\documentclass[11pt]{article}
\usepackage{amsmath,amsfonts, stix}
\usepackage{xcolor}
\numberwithin{equation}{section}
\newtheorem{theorem}{Theorem}[section]
\newtheorem{lemma}[theorem]{Lemma}
\newtheorem{proposition}[theorem]{Proposition}
\newtheorem{corollary}[theorem]{Corollary}
\newtheorem{definition}[theorem]{Definition}
\newtheorem{remark}[theorem]{Remark}

\newcommand{\footremember}[2]{%
    \footnote{#2}
    \newcounter{#1}
    \setcounter{#1}{\value{footnote}}%
}
\newcommand{\qed}{\nolinebreak\hfill\vbox{\hrule\hbox{\vrule\kern3pt\vbox{\kern6pt}\kern3pt\vrule}\hrule}}
\newenvironment{pf}{{\noindent\bf Proof.}}{\qed\newline}

\begin{document}

\title{A two phase boundary obstacle-type problem for the bi-Laplacian}

\author{Donatella Danielli\footremember{DD}{School of Mathematical and Statistical Sciences, Arizona State University. Email: ddanielli@asu.edu} \&  Alaa Haj Ali\footremember{AHA}{Department of Mathematics, Purdue University. Email: hajalia@purdue.edu}
}
\date{}
\maketitle
\begin{abstract}
In this paper we are concerned with a two phase boundary obstacle-type problem for the bi-Laplace operator in the upper unit ball.
The problem arises in connection with unilateral phenomena for flat elastic plates. It can also be seen as an extension problem to an obstacle problem for the fractional Laplacian $(-\Delta )^{3/2}$, as first observed in \cite{Y}. We establish the well-posedness and the optimal regularity of the solution, and we study the structure of the free boundary. Our proofs are based on monotonicity formulas of Almgren- and Monneau-type.
\end{abstract}

\section{Introduction}
In this paper we consider the following two-phase boundary obstacle problem

\begin{equation}\label{prb}
\left\{
\begin{array}{ll}
\Delta^2 u =0 & \text { in } B_1^+\\
u=g & \text { on } (\partial B_1)^+ \\
u_y=0 & \text { on } B_1'\\
(\Delta u)_{y}=\lambda_- (u^-)^{p-1}-\lambda_+ (u^+)^{p-1}  & \text { on }  B_1'.\\
\end{array}
\right.
\end{equation}
Here  $B_1^+ =\{z=(x,y) \in B_1 ; \ \ y>0\}$, with $B_1$ being the unit ball in $\mathbb{R}^{n+1}$ centered at the origin, $ (\partial B_1)^+= \partial B_1 \cap \{y > 0\}$, $B_1'=B_1 \cap \{y=0\}$, $\lambda_-, \lambda_+$ are positive constants, and $p>1$. For $q:=\max\{2,p\}$, we let $g\in W^{2,q}(B_1^{+})$, with $g_y=0$ on $B_1'$.

In the two dimensional case, this problem arises from Kirchhoff's phenomena of a linearized plate bending, where an isotropic thin elastic plate $\Omega$ is clamped at its edges along a portion $\Gamma_1$ of the boundary $\partial \Omega$ and satisfies a unilateral restriction on the remaining portion $\Gamma_2:=\partial \Omega \setminus \Gamma_1$. In such a model, the vertical deflection $u$ is subject to the boundary conditions  $u=g, \ \ u_{\nu}=0$ on $\Gamma_1$ (where $\nu$ is the outer unit normal vector to $\partial \Gamma_1$).  This says that the plate is clamped on a portion of its boundary along a curve $g(x)$, and it remains horizontal there. The vertical deflection of the plate is caused by a unilateral displacement applied to the portion $\Gamma_2$ of the boundary. When the displacement is small, the vertical deflection is governed by the bi-Laplace operator. For more details, we refer the reader to \cite[Chapter 2]{LL}.

This problem is equivalent to a non-local one in a lower dimensional Euclidean space, as first observed in \cite{Y}. More precisely, Yang \cite{Y} generalized the characterization of the fractional Laplacian $(-\Delta)^s$   as a Dirichlet-to-Neumann map, established by Caffarelli and Silvestre \cite{CS} for $0<s<1$,  to all positive non-integer $s$. In particular, for $s=3/2$, we recall the following result.

\begin{theorem}[{\cite[Theorem 2.1]{Y}}]\label{extension}
Suppose $u\in W^{2,2} (\mathbb{R}_{+}^{n+1})$ is a solution (in a suitable weak sense) to
\begin{equation*}
\left\{
\begin{array}{ll}
 \Delta^2u(x,y)=0 & \text { in } \mathbb{R}^n \times \mathbb{R}_{+}\\
u(x,0)=u_0(x) & \text { for all } x \in \mathbb{R}^n\\
u_{y}(x,0)=0 & \text { for all } x \in \mathbb{R}^n\\
\end{array}
\right.
\end{equation*}
for some $u_0$ in $H^{{3}/{2}}(\mathbb{R}^n)$.  Then
\begin{equation}\label{D_to_N}
(-\Delta)^{3/2} u_0(x)=C_{n} \frac{\partial}{\partial y} \Delta u(x,0).
\end{equation}
\end{theorem}

Given the above result, we see that our problem (\ref{prb}) can be seen as a local version (thanks to   the bi-harmonic extension) of a  two-phase obstacle-type problem, associated with the fractional operator $(-\Delta )^{3/2}$,
\begin{equation}\label{prb_fractional}
(-\Delta)^{3/2} u =\lambda_- (u^-)^{p-1}-\lambda_+ (u^+)^{p-1}  \ \  \text { in }  \mathbb{R}^n, \ \ \ \ u \to 0 \text { as }|x| \to \infty
\end{equation}
in the lower dimensional Euclidean space $\mathbb{R}^n$. It should be noted that in our work, however, we do not exploit this connection, and we deal exclusively with the local formulation (\ref{prb}).


The  problem (\ref{prb}) has been considered in \cite{FF}  in the case where the non-homogeneous boundary condition  is given by $(-\Delta u)_y= h(x) u(\cdot,0)$, with $h(x)$ being a $C^1$-function on $B_1'$. The authors prove a monotonicity formula of Almgren's type and classify the possible blow-up limit profiles. Then, they use their results to establish the strong unique continuation property for the fractional equation $(-\Delta)^{3/2} u(x)= a(x) u(x)$  for some $C^1$ function $a(\cdot)$.

In the present work, both the formulation of  problem (\ref{prb}) and our objectives are different. The non-homogeneous thin boundary condition $(\Delta u)_{y}=\lambda_- (u^-)^{p-1}-\lambda_+ (u^+)^{p-1}   \text { on }  B_1'$ has a right-hand-side which is non-differentiable in $u$ when $1<p\leq 2$ and non-linear in $u$ when $p>2$. Moreover, one of our main goals is to study the regularity of the solution of (\ref{prb}) for all $p>1$, and to prove that the regularity we establish is optimal in the cases where $p$ is an integer. Moreover, in the cases when $p=2$ and $p\geq 3$, we study the structural properties of the free boundary $\left(\partial \{u>0\} \cup \partial \{u<0\}\right) \cap B_1'$.

There have been several results in the literature concerning free boundary problems associated to the bi-harmonic operator. To mention a few, in \cite{DKV1} the authors study a singular perturbation problem which can be understood as the bi-harmonic counterpart of the classical combustion problem. Then, in \cite{DKV2} the same authors study the regularity of the minimizer of the Bernoulli-type functional for the bi-Laplacian: $\int_{\Omega}( |\Delta w|^2 + \chi_{w>0})dx $.
They also carry out a structural analysis of the free boundary of such a minimizer. In addition, there have been  some results related to the thick obstacle problem for the bi-Laplacian. In \cite{F1} the author establishes the $W^{3,2}_{loc}$ regularity of a solution, and in \cite{F2} and \cite{CF} the authors provide different approaches to prove the $C^{1,\alpha}_{loc}$ regularity of the solution (for $0<\alpha<1$). Finally, in the recent paper \cite{A}, the author  establishes the $C^{1,\alpha}$-regularity of the free boundary, and then improves the regularity of the solution to $C^{3,\alpha}_{loc}$, under the assumption that the solution is almost one-dimensional and its positivity set $\{u>0\}$ is a non-tangentially accessible (NTA) domain.

In the last decade, there has been a surge of interest in lower dimensional obstacle problems for second order uniformly elliptic operators. These extensive studies were motivated by the discovery of several powerful monotonicity formulas which provided the necessary tools to establish the optimal regularity of a solution, as well as the regularity and the structure of the free boundary. We refer the reader to \cite{ACS}, \cite{AC}, \cite{CSS}, \cite{CS}, \cite{GP}, the survey \cite{DS} and the book \cite{PSU}, as well as the references therein. On the other hand, much less work has been done on boundary obstacle-type problems associated with the bi-harmonic operator. In the 1980s, B. Schild studied the regularity of the solution of the thin obstacle problem for the bi-Laplacian. This problem can be seen as a limit case of our problem (in the model case $p=2$) when one of the constants $ \lambda_{\pm}$ vanishes and the other goes to $\infty$. In \cite{S1}, the author establishes the $H_{loc}^{2,\infty} \cap H_{loc}^{3/2}$ regularity of a solution, whereas in \cite{S2} he  proves the $H^{3/2}$ and $C^{2,1/2}$ optimal regularity of the solution for the two-dimensional case. In fact, the author obtains such results for the general polyharmonic operator $(-\Delta)^m$ for all integers $2 \leq m \leq \frac{n+2}{2}$.
However, many important questions related to the structure of the free boundary and the characterization of the free boundary points remained open. As far as we are aware, our work is the first instance of a comprehensive study of the free boundary in a boundary obstacle problem associated with a fourth-order operator.

The structure of the paper is as follows. In Section \ref{well_posedness}  we introduce a variational formulation of the problem to  prove the existence  and uniqueness of a minimizer for the corresponding functional. The relevant functional is similar to a thin and higher order version of the classical Alt-Phillips functional first introduced by Alt and Phillips in \cite{AlP} and \cite{P}. We then establish the weak differentiability of the Laplacian of the minimizer, and we obtain some sub-harmonicity properties for the minimizer and for its Laplacian.
In Section \ref{optimal_regularity} we establish the $C^{p+1,\alpha}$-regularity of a solution ($0<\alpha<1$) for all $p>1$, and then we show that the established regularity is optimal in the cases where $p$ is an integer.  Next, we turn our attention to obtaining all the necessary tools to study the structure  of the free boundary. In the remaining sections of the paper, for technical reasons the range of the values of $p$ is restricted  to $p=2$ or $p \geq 3$. In Section \ref{Almgren's_monotonicity} we prove an Almgren-type monotonicity formula and we deduce the growth rate of the solution away from the free boundary. In Section \ref{blow_up_analysis} we carry the analysis of the blow-up sequences around a free boundary point. In Section \ref{Monneau_monotonicity} we obtain a Monneau-type monotonicity formula and we deduce the nondegeneracy of the solution near so-called singular free boundary points. Finally, in Section \ref{structure_singular_set} we collect all the results previously obtained  to determine the structure of the singular set.

\section { Existence and Uniqueness of a minimizer} \label{well_posedness}

We begin by introducing the notion  of a weak solution for the problem (\ref{prb}). We recall $q:=\max\{2,p\}$.
\begin{definition}
We say that $u \in W^{2,q}(B_1^+)$ is a weak solution of (\ref{prb})  if
\begin{equation}\label{weaksol}
 \int_{B_1^+} \Delta u \Delta \phi dx dy =\int_{B_1'}\left(\lambda_- (u^-)^{p-1}  -\lambda_+ (u^+)^{p-1}\right)\phi dx
 \end{equation}
  for all functions $\phi \in C^{\infty} (B_1^+)$ such that  $\phi=0 \text { on }  (\partial B_1)^+$ and $\phi_y=0 \text { on }  B_1'$.
\end{definition}
To prove the existence of a weak solution to (\ref{prb}), we follow a variational approach by minimizing the functional
\begin{equation}\label{fnct}
J[w]:=\int_{B_1^+}  (\Delta w)^2 dx dy +  \frac{2}{p} \int_{B_1'} \lambda_- (w^-)^p+\lambda_+ (w^+)^p dx
\end{equation}
over the admissible set
\begin{equation}\label{admissible}
\mathcal{A}:= \left\{ w \in W^{2,q} (B_1^+), \ \ w=g \text { on } (\partial B_1)^+, \text { and } \ \ w_y=0 \text { on } B_1'\right\}.
\end{equation}
This is a free boundary problem with two phases. We will denote the free boundary with
\begin{equation}\label{free_bndry}
\Gamma=\Gamma_- \cup \Gamma_+ \text{, where }
\Gamma_-=\partial_x (\{u<0\} \cap B_1')
\text{ and }\Gamma_+=\partial_x (\{u>0\} \cap B_1').
\end{equation}

Our starting point is the following Poincar\'e type inequality on a half ball. This result was proven in \cite[Lemma 2.3]{FF}, but for the reader's  convenience we include the proof here.

\begin{lemma}\label{poincareineq*}
For $w \in W^{1,2}(B_r)$ we have
\begin{equation}\label{poincareineq}
\frac{n}{r^2} \int_{B_r^+} w^2 dx dy \leq \frac{1}{r}\int_{(\partial B_r)^+} w^2 dS + \int_{B_r^+} |\nabla w|^2 dx dy.
\end{equation}
\end{lemma}
\begin{pf}
Letting $z:=(x,y)$ we have
\begin{equation*}
\begin{split}
(n+1) \int_{B_r^+} w^2 dx dy =& \int_{B_r^+}   \nabla \cdot (w^2 z)-2\int_{B_r^+} w \nabla w \cdot z=r \int_{(\partial B_r)^+} w^2 -2\int_{B_r^+} w \nabla w \cdot z\\
& \leq r \int_{(\partial B_r)^+} w^2+ \int_{B_r^+} w^2 dx dy + r^2\int_{B_r^+} |\nabla w|^2 dx dy,\\
\end{split}
\end{equation*}
which gives
\begin{equation*}
n \int_{B_r^+} w^2 dx dy\leq r \int_{(\partial B_r)^+} w^2 dS+ r^2\int_{B_r^+} |\nabla w|^2 dx dy.
\end{equation*}
Dividing by $r^2$ we obtain the desired inequality.
\end{pf}

We proceed to establish existence and uniqueness of the minimizer.
\begin{theorem}\label{exist_min}
There exists a unique minimizer $u$ of the functional  (\ref{fnct}) over the admissible set $\mathcal{A}$.
\end{theorem}
\begin{pf}
First we note that the admissible set $\mathcal{A} \neq \emptyset$ since $g \in \mathcal{A}$. Next, let $\{u_l\}$  be a minimizing sequence of $J[\cdot]$ and denote $v_l:=\Delta u_l$ for all $l \geq 1$. Clearly, $\{v_l\}$ is bounded in $\mathcal{L}^2(B_1^+)$. Moreover, integrating by parts and applying  (\ref{poincareineq}), we have
\begin{equation*}
\begin{split}
&\int_{B_1^+} |\nabla (u_l-u_1)|^2 dx dy =-\int_{B_1^+} (u_l-u_1) (v_l-v_1)dx dy \leq ||u_l-u_1||_{\mathcal{L}^2(B_1^+)} ||v_l-v_1||_{\mathcal{L}^2(B_1^+)}\\
&\leq C||\nabla (u_l- u_1)||_{\mathcal{L}^2(B_1^+)} ||v_l-v_1||_{\mathcal{L}^2(B_1^+)},\\
\end{split}
\end{equation*}
and thus
\begin{equation*}
||\nabla (u_l- u_1)||_{\mathcal{L}^2(B_1^+)} \leq C ||v_l-v_1||_{\mathcal{L}^2(B_1^+)}.
\end{equation*}
Also,
$$||u_l- u_1||_{\mathcal{L}^2(B_1^+)} \leq C||\nabla (u_l- u_1)||_{\mathcal{L}^2(B_1^+)}\leq C^2 ||v_l-v_1||_{\mathcal{L}^2(B_1^+)}.$$
Moreover, by classical elliptic regularity we have
\begin{equation*}
||u_l-u_1||_{W^{2,2}(B_1^+)} \leq C\left(||u_l-u_1||_{\mathcal{L}^2(B_1^+)}+||v_l-v_1||_{\mathcal{L}^2(B_1^+)} \right)
\end{equation*}
for some constant $C$ depending only on the dimension.
Therefore, the sequence $\{u_l\}$ is bounded in $W^{2,2}(B_1^+)$.
By standard compactness arguments, for a subsequence still denoted by $\{u_l\}$ we have
\begin{equation*}
u_l \to u \text{ weakly in } W^{2,2}(B_1^+) \text { and strongly in } \mathcal {L}^{2}(B_1^+).
\end{equation*}
Passing to a subsequence (without changing its name for simplicity), we may assume that $u_l$ converges a.e. in $B_1^+$. Since the admissible set $\mathcal{A}$ is a closed subset of $W^{2,2}(B_1^+)$,  we conclude that $u\in\mathcal{A}$. From the weak convergence  we infer
\begin{equation}\label{term1}
\int_{B_1^+}(\Delta u)^2 dx dy \leq\liminf_{l\to\infty}\int_{B_1^+}(\Delta u_l)^2 dx dy.
\end{equation}
Now we claim that we also have
\begin{equation}\label{term2}
\int_{B_1'} \lambda_- (u^-)^p+\lambda_+ (u^+)^p dx \leq \liminf_{l\to\infty} \int_{B_1'} \lambda_- (u_l^-)^p+\lambda_+ (u_l^+)^p dx.
\end{equation}
To see this, we note that since the half unit ball is a Lipschitz domain, the trace operator
$T:W^{2,2}(B_1^+) \rightarrow \mathcal{L}^2(\partial B_1^{+})$ is a compact bounded linear operator. Therefore, $u_l \rightarrow u$ strongly in $\mathcal{L}^2(B_1')$ and, for a subsequence,  the convergence  is a.e. in $B_1'$.
Clearly, this also implies that $\lambda_- (u_l^-)^p+\lambda_+ (u_l^+)^p$ converges a.e. to $\lambda_- (u^-)^p+\lambda_+ (u^+)^p$ in $B_1'$. Finally, an application of Fatou's lemma finishes the proof of claim (\ref{term2}). Combining (\ref{term1})  with (\ref{term2}) we obtain $$J [u] \leq \lim \inf J[u_l], $$
which implies that $u$ is a minimizer of the functional $J[\cdot]$.
The uniqueness of the minimizer follows from the convexity of our functional.
\end{pf}

\begin{theorem}\label{exist_sol}
The minimizer of the functional (\ref{fnct}) over the set $\mathcal{A}$ is a weak solution of (\ref{prb}).
\end{theorem}

\begin{pf}
Fix a function $\phi \in W^{2,2}(B_1^+)$ such that $\phi=0 \text { on } (\partial B_1)^+$ and $\phi_y=0 \text { on }  B_1'$. Writing $\phi=\phi^{+}-\phi^{-}$ and noting that  both $\phi^{+}$ and $\phi^{-}$ are valid test functions, without loss of generality we can assume  $\phi~\geq~0$. Clearly, $u+\tau \phi$ is an admissible function for all values of $\tau$, and by the minimality of $u$ we have $J[u+\tau \phi] \geq J[u]$. Then, for all $\tau>0$, we have
\begin{equation}\label{+per}
\begin{split}
 0 \leq\frac{J[u + \tau \phi]-J[u]}{\tau}&=\int_{B_1^+}  \left((\Delta \phi)^2 \tau+ 2 \Delta u \Delta \phi \right)\\
 & +\frac{2}{p} \int_{B_1'}\frac
 {\left(\lambda_- ((u+\tau \phi)^-)^p+\lambda_+ (((u+\tau \phi))^+)^p\right)-\left(\lambda_- (u^-)^p+\lambda_+ ((u)^+)^p\right)}{\tau }.\\
 \end{split}
\end{equation}
By convexity we see that the integrand in the last line is monotone in $\tau$ and converges
to $$p\left( \lambda_+ (u^+)^{p-1}-\lambda_- (u^-)^{p-1} \right) \phi\qquad \mbox{as } \tau\searrow 0.$$
Therefore, letting $\tau \searrow 0$ and applying the monotone convergence theorem we get
\begin{equation*}
\int_{B_1^+} \Delta u \Delta \phi + \int_{B_1'}  \lambda_+ (u^+)^{p-1}-\lambda_- (u^-)^{p-1} \phi \geq  0.
\end{equation*}
Similarly, an opposite inequality  to (\ref{+per}) holds for $\tau<0$. Letting $\tau \nearrow 0$, we obtain the inequality
\begin{equation*}
\int_{B_1^+} \Delta u \Delta \phi + \int_{B_1'}  \lambda_+ (u^+)^{p-1}-\lambda_- (u^-)^{p-1}\phi  \leq  0,
\end{equation*}
and this finishes the proof of the theorem.
\end{pf}

From now on, we will denote
\begin{equation}\label{v_definition}
v:=\Delta u
\end{equation}
The following theorem shows the weak differentiability of $v$.

\begin{theorem}\label{gradv}
The function $v$  is weakly differentiable, and its weak derivatives are in $\mathcal{L}(B^+_{\rho})$ for all balls $B_{\rho} \subset \subset B_1$.
\end{theorem}
\begin{pf}
Let $\eta$ be a smooth cut-off function in $B_1$ such that $\eta \equiv 1$ on $B_{1/2}$. We wish to show that $\eta v$ admits a weak derivative in $B_{1/2}$. To this end, we let $w=\eta v$ and we note that $w$ satisfies
\begin{equation}\label{w_est}
\begin{split}
&\int_{B_1^+} w \Delta \phi dx dy=\int_{B_1^+} v\eta \Delta \phi dx dy= -\int_{B_1^+}\left(2 \nabla \eta \cdot \nabla \phi + \phi \Delta \eta \right) v+\int_{B_1^+} v \Delta (\eta \phi)\\
&=-\int_{B_1^+} 2 v \nabla \eta \cdot \nabla \phi -\int_{B_1^+} \Delta \eta v \phi+\int_{B_1'}\left(\lambda_- (u^-)^{p-1}  -\lambda_+ (u^+)^{p-1}\right)\eta \phi dx
 \end{split}
 \end{equation}
 for all functions $\phi \in C^{\infty} (B_1^+)$ such that  $\phi=0 \text { on }  (\partial B_1)^+$ and $\phi_y=0 \text { on }  B_1'$. That is,  $w$ is a solution of the problem
\begin{equation}\label{prb1}
\left\{
\begin{array}{ll}
\int_{B_1^+} w \Delta \phi dx dy=-\int_{B_1^+} 2v \nabla \eta \cdot \nabla \phi -\int_{B_1^+} \Delta \eta v \phi+\int_{B_1'}\left(\lambda_- (u^-)^{p-1}  -\lambda_+ (u^+)^{p-1}\right)\phi dx &\text { in } B_1^+\\
w=0 & \text { on } (\partial B_1)^+\\
\end{array}
\right.
\end{equation}
Moreover, a solution to (\ref{prb1}) is unique. To see this, suppose $w_1$ and $w_2$ are two solutions of (\ref{prb1}) and let $w_0:=w_2-w_1$. Then, $w_0$ satisfies
\begin{equation}\label{w0_eq1}
\int_{B_1^+} w_0 \Delta \phi dx dy=0
\end{equation}
for all functions $\phi \in C^{\infty} (B_1^+)$ such that  $\phi=0 \text { on }  (\partial B_1)^+$ and $\phi_y=0 \text { on }  B_1'$.
Now we show that we also have
\begin{equation}\label{w0_eq2}
\int_{B_1^+} w_0  \psi dx dy=0
\end{equation}
for all functions $\psi \in C^{\infty} (B_1^+)$ such that  $\psi=0 \text { on }  (\partial B_1)^+$, and this would immediately imply that $w_0=0$ in $B_1^+$. To this end, we observe that for  such a function $\psi$,  the problem
 \begin{equation}\label{prb2}
\left\{
\begin{array}{ll}
\Delta \phi=\psi &\text { in } B_1^+\\
\phi=0 & \text { on } (\partial B_1)^+\\
\phi_y=0 & \text { on } \partial B_1'\\
\end{array}
\right.
\end{equation}
admits an unique solution $\phi$. This in turn implies that $\psi$ satisfies (\ref{w0_eq2}), and shows that solutions to (\ref{prb1}) are unique.

On the other hand, for $\eta$, $u$ and $v$ fixed as above, the  right-hand-side of the PDE in (\ref{prb1}) is a bounded linear functional on $H_0^{1}(B_1)$.
Therefore, by the Lax-Milgram theorem, there exists a function $\Psi$ in $H_0^1(B_1)$ such that
\begin{equation*}
\int_{B^+_1} \nabla \Psi \nabla \phi dx dy=\int_{B^+_1} 2v \nabla \eta \cdot \nabla \phi +\int_{B^+_1} \Delta \eta v \phi+\int_{B_1'}\left(\lambda_+ (u^+)^{p-1}-\lambda_- (u^-)^{p-1}\right)\phi dx
\end{equation*}
for all $\phi \in C_0^{\infty}(B_1)$.
Now integrating by parts in the left-hand-side of the equation above, we see that $\Psi$ is a solution of the problem (\ref{prb1}).
Finally, by the uniqueness of a solution of the problem (\ref{prb1}), we conclude that $\Psi=\eta w$ in $B^+_1$. Thus, $\Psi=v$ in $B^+_{1/2}$ and the result follows.
\end{pf}

\begin{proposition}
The minimizer $u$ satisfies
\begin{equation}\label{loc_weaksol}
\int_{B_{\rho}^+} \nabla v \cdot \nabla \phi + \int_{B_{\rho}{'}}\left(\lambda_- (u^-)^{p-1}-\lambda_+ (u^+)^{p-1}\right)\phi dx=0
\end{equation}
for all test functions $\phi \in C^{\infty}(B^+_{\rho})$ and satisfying $\phi=0$ on $(\partial B_{\rho})^+$.
In particular, $v$ is harmonic in $B_1^{+}$.
\end{proposition}
\begin{pf}
Thanks to Theorem \ref{gradv}, an integration by parts of (\ref{weaksol}) implies that (\ref{loc_weaksol}) holds for all test functions $\phi$ which satisfy the additional restriction $\phi_y=0$ on $B_1'$.
Now, for a general test function $\phi$, consider the sequence of test functions $$\phi_{\epsilon}(x,y):=\phi(x,y)-\phi_y(x,y)y\eta\left(\frac{y}{\epsilon}\right), $$
where $\eta$ is a smooth function on $\mathbb{R}$, compactly supported in $(-2,2)$ and such that $\eta \equiv 1$ in $(-1,1)$. It is easy to see that  $\phi_{\epsilon} \to \phi$ and $\nabla \phi_{\epsilon} \to \nabla \phi$ uniformly as $\epsilon \to 0$. Moreover, $\phi_{\epsilon}$ is a test function which satisfies the additional restriction $\phi_y=0$ on $B_1'$. Thus,
$$\int_{B_1^+} \nabla v \cdot \nabla \phi_{\epsilon}+ \int_{B_1'}\left(\lambda_- (u^-)^{p-1}-\lambda_+ (u^+)^{p-1}\right)\phi_{\epsilon} dx=0$$
and the result follows by passing to the limit as $\epsilon \to 0$.
\end{pf}

\begin{proposition}\label{v_min}
For every $0<\rho <1$, the function $v$ is a minimizer of the functional
\begin{equation}\label{v_fnct}
I[w]:=\frac{1}{2} \int_{B_{\rho}^+} |\nabla w|^2 +\int_{B_{\rho}{'}}\left(\lambda_- (u^-)^{p-1}-\lambda_+ (u^+)^{p-1}\right)w dx
\end{equation}
over the set $$\mathcal{A}_{\rho}:= \left\{w\in W^{1,2}(B_{\rho}^+) \ \ , \ \ w=v \text { on } (\partial B_{\rho})^+\right\}.$$
\end{proposition}
\begin{pf}
For $w \in \mathcal{A}_{\rho}$ we have
\begin{equation*}
I[w]-I[v]=\frac{1}{2} \int_{B_{\rho}^+} |\nabla (w-v)|^2 +\int_{B_{\rho}^+} \nabla v \cdot \nabla (w-v)+\int_{B_{\rho}{'}}\left(\lambda_- (u^-)^{p-1}-\lambda_+ (u^+)^{p-1}\right)(w-v) dx\geq 0,
\end{equation*}
where the last inequality follows from (\ref{loc_weaksol}).
\end{pf}

 In the remaining of the paper, we extend $u$ to the entire unit ball $B_1$ by even symmetry.
\begin{lemma}\label{sub_super_sol}
Let $u$ be the minimizer of (\ref{fnct}). Then,
$u^{\pm}$ and $v^{\pm}$ are subharmonic in $B_1$.
\end{lemma}

\begin{pf}
Fix a non-negative test function $\phi$ such that $\phi=0$ on $(\partial B_1)^{+}$ and $\phi_y=0$ on $B_1'$. For small $\epsilon>0$, let $u_{\epsilon} = (u-\epsilon \phi)^+-u^-$. Clearly,  $u_{\epsilon}$ is an admissible function. Also, we see that $u_{\epsilon}^+ \leq u^+$ and $u_{\epsilon}^-=u^-$. This, along with the minimality of $u$, yields
\begin{equation*}
\int_{B^+_1 \cap \{u>0\}} (\Delta u)^2 dx dy \leq  \int_{B^+_1\cap\{u>\epsilon \phi\}} (\Delta u-\epsilon \Delta \phi)^2 dx dy \leq  \int_{B^+_1\cap\{u>0\}} (\Delta u-\epsilon \Delta \phi)^2 dx dy,
\end{equation*}
which gives
\begin{equation}\label{u+varineq}
\int_{B^+_1} \Delta u^+ \Delta \phi \leq 0.\\
\end{equation}
Now, for a non-negative smooth  function $\psi$ such that $\psi=0$ on $(\partial B_1)^+$, let $\phi$ be the unique solution to
\begin{equation}\label{prb2*}
\left\{
\begin{array}{ll}
\Delta \phi=\psi &\text { in } B_1^+\\
\phi=0 & \text { on } (\partial B_1)^+\\
\phi_y=0 & \text { on } \partial B_1'\\
\end{array}
\right.
\end{equation}
By the maximum principle we have $\phi \leq 0$ on $B_1^+$. Therefore, $-\phi$ is admissible for  (\ref{u+varineq}), which immediately implies
\begin{equation*}
\int_{B^+_1} \Delta u^+ \psi \geq 0.\\
\end{equation*}
Integrating by parts we obtain
\begin{equation}\label{u+sub}
\int_{B^+_1} \nabla u^+ \cdot \nabla \psi \leq 0,\\
\end{equation}
that is $u^+$ is subharmonic in $B_1$.
Arguing in a similar way we can also show that $u^-$ is subharmonic in $B_1^+$.

Next, for a fixed $0<\rho<1$, take a non-negative test function $\phi$ such that $\phi=0$ on $(\partial B_{\rho})^{+}$, and for small $\epsilon>0$ let $v_{\epsilon} = (v)^+-(v+\epsilon \phi)^-$. Clearly  $v_{\epsilon}\in \mathcal{A}_{\rho}$,  $v_{\epsilon}^- \leq v^-$ and $v_{\epsilon}^+=v^+$. Thus,  by Proposition~\ref{v_min} we have
\begin{equation*}
\begin{split}
\frac{1}{2} \int_{B_{\rho}^+\cap \{v<0\}} |\nabla v|^2 &+\int_{B_{\rho}{'}\cap \{v<0\}} \lambda_- (u^-)^{p-1} v^- dx\\
 &\leq \frac{1}{2} \int_{B_{\rho}^+\cap \{v<-\epsilon \phi\}} |\nabla v+\epsilon \nabla \phi|^2 +\int_{B_{\rho}{'} \cap \{v<-\epsilon \phi\}} \lambda_- (u^-)^{p-1} (v+\epsilon \phi)^- dx\\
&\leq \frac{1}{2} \int_{B_{\rho}^+\cap \{v<0\}} |\nabla v+\epsilon \nabla \phi|^2 +\int_{B_{\rho}{'} \cap \{v<0\}} \lambda_- (u^-)^{p-1} v^- dx\\
\end{split}
\end{equation*}
for all $\epsilon>0$. Hence,
\begin{equation}\label{v_sub}
\int_{B_{\rho}^+} \nabla v^- \cdot \nabla \phi \leq 0
\end{equation}
and $v^-$ is subharmonic in $B_{\rho}$. We can show that the same result holds for $v^+$ in a similar way.
\end{pf}

\begin{remark}
We have $v=0$ on the free boundary $\Gamma$ defined in (\ref{free_bndry}).
\end{remark}

\begin{pf}
This is an immediate consequence of the definition of a weak solution (\ref{weaksol}) and Lemma \ref{sub_super_sol}, which implies that $v^{\pm}=\Delta u^{\pm}$.
\end{pf}

\section{Regularity of the minimizer}\label{optimal_regularity}
We recall that we have extended our solution $u$ to the entire unit ball $B_1$ by even symmetry. We first prove the following local boundedness result.

\begin{lemma}\label{localbnd}
Let $u$ be the minimizer to (\ref{fnct}). Then both  $u$ and its Laplacian $v$ are in $\mathcal{L}^{\infty}_{loc}(B_1)$. 
\end{lemma}
\begin{pf}
Let $B_R \subset \subset B_1$. By Lemma \ref{sub_super_sol} and the $L^{\infty}-L^2$ estimate for subharmonic functions, we have
\begin{equation}\label{u_loc_bnd}
||u||_{\mathcal{L}^{\infty}(B_R)} \leq C(R) ||u||_{\mathcal{L}^{2}(B_1)}
\end{equation}
and
\begin{equation}\label{v_loc_bnd}
||v||_{\mathcal{L}^{\infty}(B_R)} \leq C(R) ||v||_{\mathcal{L}^{2}(B_1)},
\end{equation}
where $C(R)$ is a constant which depends only on $R$ and the dimension $n$.\end{pf}

In the remainder of the paper we denote
\begin{equation}\label{u_R_bnd}
m_1(R):=||u||_{\mathcal{L}^{\infty}(B_R)}
\end{equation}
and
\begin{equation}\label{v_R_bnd}
m_2(R):=||v||_{\mathcal{L}^{\infty}(B_R)}.
\end{equation}

Next, we obtain the following energy estimate.

\begin{lemma}\label{enrg_est}
Let $u$ be the minimizer of (\ref{fnct}). Then, for any $B_{2r} \subset \subset B_1$ we have the estimates

\begin{equation}\label{u_erg_est}
\int_{B_{r}} |\nabla u|^2 \leq \frac{C}{r^2} \int_{B_{2r}} u^2,
\end{equation}
where the constant $C$ depends only on the dimension $n$.
\end{lemma}
\begin{pf} The energy inequality (\ref{u_erg_est}) holds given that $u^{\pm}$ are subharmonic in $B_1$.
\end{pf}

The following lemma gives a regularity result for our solution $u$. Then, in Theorem \ref{optimal_reg}, we will prove that this regularity result is optimal in the case where $p$ is an integer.

\begin{lemma}\label{part_reg}

Let $u$ be the minimizer of (\ref{fnct}). When $p$ is an integer,
 $u \in C_{loc}^{p+1, \alpha}(B_1^+\cup B_1')$ for all $\alpha<1$. Moreover,
 \begin{equation}\label{u_part_reg_norm1}
||u||_{C^{p+1,\alpha}(B_{\rho}^+\cup B_{\rho}^{'}) } \leq C(\alpha,\rho) \left(||u||_{\mathcal{L}^{\infty}(B_{\rho})}+||v||_{\mathcal{L}^{\infty}(B_{\rho})}\right).
\end{equation}
When $p$ is not an integer, $u \in C_{loc}^{\lfloor p+1 \rfloor, \alpha}(B_1^+)\cup B_1'$  for all $\alpha<p-1-\lfloor p-1 \rfloor$ and
\begin{equation}\label{u_part_reg_norm2}
||u||_{C^{p+1,\alpha}(B_{\rho}^+\cup B_{\rho}^{'})} \leq C(\alpha,\rho) \left(||u||_{\mathcal{L}^{\infty}(B_{\rho})}+||v||_{\mathcal{L}^{\infty}(B_{\rho})}\right).
\end{equation}
\end{lemma}
\begin{pf}
As a first step, we show that $u$ is locally Lipschitz continuous in $B_1$. Fix a point $x$ on $B_1'$ and a ball $B_R$ centered at $x$. For $0<r<R$, let $w$ be the harmonic replacement of $u$ in $B_r$. We then have
\begin{equation}\label{u_w_est1}
\int_{B_r} |\nabla u-\nabla w|^2 dx dy= \int_{B_r} \nabla u \cdot (\nabla u-\nabla w) dx dy=-\int_{B_r} v (u-w) dx dy,
\end{equation}
where the first equality follows from $w$ being harmonic, and the second one from an integration by parts. Now, an application of the maximum principle to the harmonic function $w$ on $B_r$ yields
$$ |w| \leq ||u|| _{\mathcal{L}^{\infty}(B_R)}=m_1(R).$$
This, along with (\ref{u_w_est1}), gives
\begin{equation}\label{approx}
\int_{B_r} |\nabla u-\nabla w|^2 dx dy \leq C(R) r^{n+1}.
\end{equation}
Here, $C(R)$ depends only on $m_1(R)$, $m_2(R)$ and the dimension $n$.

Next, fix $0<r<R$. For all $0<\rho<r$ we have
\begin{equation*}
\begin{split}
 \left(\int_{B_{\rho}(x)} |\nabla u| ^2\right)^{\frac{1}{2}}&\leq C\left(\int_{B_{r}(x)} |\nabla (u-w)| ^2\right)^{\frac{1}{2}} + \left(\int_{B_{\rho}(x)} |\nabla w|^2\right)^{\frac{1}{2}}\\
&\leq C(R) r^{\frac{n+1}{2}} +\left(\frac{\rho}{r}\right)^{\frac{n+1}{2}}\left(\int_{B_{r}(x)} |\nabla w|^2\right)^{\frac{1}{2}}\leq C(R) r^{\frac{n+1}{2}} + \left(\frac{\rho}{r}\right)^{\frac{n+1}{2}}\left(\int_{B_{r}(x)} |\nabla u|^2\right)^{\frac{1}{2}}.
\end{split}
\end{equation*}
The second inequality follows from (\ref{approx}) and the harmonicity of $w$, whereas the last one is a consequence of the minimality property of the harmonic function $w$. Using the energy inequality (\ref{u_erg_est}) and arguing as in \cite[Theorem 3.1]{AP} we arrive at
\begin{equation*}
 \int_{B_{\rho}(x)} |\nabla u| ^2\leq C(R) \rho^ {n+1}
 \end{equation*}
for all $0<\rho<r$. Applying Morrey's inequality we conclude that $u$ is $C_{loc}^{1,1/2}$. In particular, $u$ is locally Lipschitz continuous in $B_1$. This in turn implies that,  in the case $1<p \leq 2$,  the function $\lambda_- (u^-)^{p-1}-\lambda_+ (u^+)^{p-1}$ is H\"older continuous of order $\alpha$  for all $\alpha<p-1$. In the case $p>2$, we conclude that $u \in C_{loc}^{p-2, \alpha}$ for all $\alpha<1$ when $p$ is an integer, and  $u \in C_{loc}^{\lfloor p-2 \rfloor, \alpha}$ for all $\alpha<p-1-\lfloor p-1 \rfloor$ when $p$ is not an integer. Since $v$ is a weak solution on $B^+_{r}$ to an oblique derivative problem, applying the regularity theory of such a problem we conclude that, when $p$ is an integer, $v$ is $C_{loc}^{p-1, \alpha}$ in $B_1^+$   up the boundary $B_1'$ for all $\alpha<1$. Moreover, when $p$ is not an integer, $v \in C_{loc}^{\lfloor p-1 \rfloor, \alpha}$ in $B_1^+$ up the boundary $B_1'$ for all $\alpha<p-1-\lfloor p-1 \rfloor$. Finally, since $u$ solves
$$\Delta u=v \text { in } B_1^+ \ \ , \ \ u_y=0 \text{ on } B_1',$$ applying the same regularity theory we reach the desired conclusion.
\end{pf}

Now we prove that the $C_{loc}^{p+1,\alpha}$ regularity is optimal. More specifically, we show that  in the case where $p$ is an integer, $u$ is not $C^{p+1,1}$ near free boundary points $z_0$ with $\nabla_{x} u (z_0) \neq 0$.

\begin{theorem}\label{optimal_reg}
In the case where $p$ is an integer, let $u$ be the minimizer of the functional (\ref{fnct}) and $z_0 \in \Gamma$ a free boundary point such that $\nabla_x u (z_0) \neq 0$. Then $u$ is not $C^{p+1,1}$ at $z_0$.
\end{theorem}

\begin{pf}
After a change of coordinates,
without loss of generality we may assume that $z_0=0$, $u_{x_1}(0)=b \neq 0$ and $u_{x_i}(0)=0$ for all $2\leq i \leq n$. We argue by contradiction. Suppose that there exists $\delta>0$ such that $u$ is in $C^{p+1,1}(B_{\delta}(0))$. For the sake of simplicity, we will assume $\delta=1$. Then,
we have
\begin{equation}\label{assump}
||D^{p+2} u||_{\mathcal{L}^{\infty} (B^+_1)}<\infty.
\end{equation}
Let $P(z)$ be the degree $(p+1)-$Taylor polynomial of $u$ and consider the family of functions
\begin{equation*}
u_{\rho} (z):= \frac{u(\rho z)-P(\rho z)}{\rho^{p+2}}
\end{equation*}
for all $\rho>0$ and  $z \in B_{1/\rho}^+$. Clearly  $u_{\rho} \in C^{p+1,1}(B_{1/\rho}^+)$ for all $\rho>0$. Also, the function $u_{\rho}$ satisfies
$$\Delta^2 u_{\rho}=-\frac{1}{\rho^{p-2}} \Delta^2 P(\rho z) \ \ \text { in }B_{1/\rho}^+.$$
We recall that $u_y(0)=v_y(0)=0$. Then $P_y=0$ and $(\Delta P)_y=0$ on $\mathbb{R}^n \times \{0\}$. We thus have
$$(u_{\rho})_y=0 \text{ on } B'_{1/\rho}$$
and
\begin{equation*}\begin{split}
(\Delta u_{\rho})_y&=\frac{v_y(\rho z)}{\rho^{p-1}}\\
&=\frac{\lambda_- (u^-(\rho z))^{p-1}-\lambda_+ (u^+(\rho z))^{p-1} }{\rho^{p-1}}\\
&=\lambda_- \left(\left(\frac{P(\rho z)}{\rho}+\rho^{p+1} (u_{\rho})\right)^-\right)^{p-1}-\lambda_+ \left(\left(\frac{P(\rho z)}{\rho}+\rho^{p+1} (u_{\rho})\right)^+\right)^{p-1}
\end{split}\end{equation*}
on $B'_{1/\rho}$. Moreover, from (\ref{assump}) we deduce the uniform bound
\begin{equation}\label{uniform_bnd}
||u_{\rho}||_{C^{p+1, 1} (B^+_{1/\rho})} \leq C(n,p) ||D^{p+2} u||_{\mathcal{L}^{\infty} (B_1^+)}.
\end{equation}
By the Arzel\`a-Ascoli theorem, there exists a subsequence, still denoted by $u_{\rho}$, and a function $w$ such that  $u^{\rho} \to w$ in the $C^{p+1}$-topology on compact subsets of $\overline{\mathbb{R}^{n+1}_+}$ as $\rho \to 0^+$. Now we have
 \begin{equation*}
\Delta^2 u_{\rho} (z)= \frac{-\Delta^2 P(\rho z)}{\rho^{p-2}} \longrightarrow Q(z):=\sum_{|\alpha|=p-2}\frac{1}{\alpha!}D^{\alpha} (\Delta^2 P)(0) z^{\alpha}
\end{equation*}
uniformly of compact sets of $\overline{\mathbb{R}^{n+1}_+}$ as $\rho \to 0^+$. Also, $\frac{P(\rho x)}{\rho} \to \nabla_x u(0) \cdot x=b x_1$. Letting $v_0=\Delta w$, we see that $v_0$ is a $C^{p-1,1}(\mathbb{R}_+^{n+1})$ solution of the problem
\begin{equation*}
\left\{
\begin{array}{ll}
\Delta v_0=Q(z) & \text { in } \mathbb{R}^{n+1}_+\\
(v_0)_y=\lambda_- (b x_1^-)^{p-1}-\lambda_+ (b x_1^+)^{p-1} & \text { on } \mathbb{R}^n \times\{0\} \\
\end{array}
\right.
\end{equation*}
The rest of the proof follows the lines of the one of \cite[Theorem 4.2]{DK} to show that such a function $v_0$ does not exist. This gives a contradiction and finishes the proof of the theorem.
\end{pf}

In the remainder of the paper, we assume that either $p=2$ or $p\geq 3$.
\section{Almgren's Type Frequency Formula}\label{Almgren's_monotonicity}

For $u, v\in W^{1,2}_{loc}(B_1)$ we introduce the \emph{Almgren's Frequency Functional}
\begin{equation}\label{Almgren's_frq}
N_0(r,u,v):=r \frac{\int_{B_r^+} (|\nabla u|^2+|\nabla v|^2) dx dy}{\int_{(\partial B_r)^+} (u^2+ v^2) dS}.
\end{equation}
For $u$  the solution of our problem (\ref{prb}), we aim to show that $N_0(r):=N_0( r, u, \Delta u)$ has a non-negative limit as $r \to 0^+$. To this end, we first show that  a suitable perturbation of $N_0(r)$  has a limit  as $r \to 0^+$. More precisely, for $u, v\in W^{1,2}_{loc}(B_1)$, we define the \emph{perturbed Almgren's Frequency Functional}
\begin{equation}\label{Almgren's_frq_per}
N(r,u,v):=r \frac{\int_{B_r^+} (|\nabla u|^2+|\nabla v|^2 + uv) dx dy+\int_{B_r{'}} \left(\lambda_- (u^-)^{p-1}  -\lambda_+ (u^+)^{p-1}\right)v dx}{\int_{(\partial B_r)^+} (u^2+ v^2) dS}.
\end{equation}
Our first step consists in proving the following result.
\begin{theorem}\label{Almgren's_per}
For $u$ solution to (\ref{prb}) with $p=2$ or $p\geq 3$, let $N(r):=N( r, u, \Delta u)$. Then,
$$\mu:=\lim_{r \to 0^{+}} N(r) \text { is finite and } \mu \geq 0.$$
\end{theorem}
The proof follows the approach used in \cite{FF}. We introduce
$$D(r):=\int_{B_r^+} (|\nabla u|^2+|\nabla v|^2 + uv) dxdy +\int_{B_r{'}} \left(\lambda_- (u^-)^{p-1}  -\lambda_+ (u^+)^{p-1}\right)v dx,$$
$$H(r):=\int_{(\partial B_r)^+} (u^2+v^2) dS$$
and
$$D_0(r):=\int_{B_r^+} (|\nabla u|^2+|\nabla v|^2 ) dxdy.$$
We will need the following two auxiliary results.
\begin{lemma}[Rellich Identity]
For $w\in W^{1,2}(B_r^+)$ we have
\begin{equation}\label{Rellich}
r \int_{(\partial B_r)^+} (|\nabla w|^2-2w_r^2)dS=(n-1) \int_{B_r^+} |\nabla w|^2 -2\int_{B_r^+} (x\cdot \nabla w)\Delta w-2\int_{B_r{'}} (x\cdot \nabla w)w_y.
\end{equation}
\end{lemma}
\begin{pf} This is a classical result, based on integration by parts and the divergence theorem.
\end{pf}
\begin{lemma}[Trace inequality]\label{l.trace}
For $u \in W^{1,2}(B_r^+)$ and $B_r \subset B_1$  we have
\begin{equation}\label{trace}
\int_{B_r{'}} u^2 dx \leq C\left( r \int_{B_r^+} |\nabla u|^2 +\int_{(\partial B_r)^+} u^2 dS  \right).
\end{equation}
\end{lemma}
\begin{pf} We refer to \cite[Lemma 2.2]{FF}.
\end{pf}

Next, we show  that
\begin{equation}\label{Nliminf}
\liminf_{r \to 0} N(r)\geq 0.
\end{equation}
To see this, we note that in view of Lemma \ref{poincareineq} we have
\begin{equation}\label{extra1}
\left|\int_{B_r^+} u v dx dy \right| \leq C \int_{B_r^+} (u^2+ v^2) dx dy \leq C\left(r^2 D_0(r)+r H(r)\right).\\
\end{equation}
Moreover, we infer from Lemma \ref{l.trace} that
\begin{equation}\label{extra2}
\left|\int_{B_r{'}} \left(\lambda_- (u^-)^{p-1}  -\lambda_+ (u^+)^{p-1}\right)v dx\right| \leq C \left(rD_0(r) + H(r)\right),
\end{equation}
where the constant $C$ depends also on $m_1(R)$ when $p>2$.
We thus have
\begin{equation}\label{N_to_N0}
N(r)\geq r\frac{D_0(r)-CrD_0(r)-CH(r)}{H(r)}\geq  (1-Cr) r\frac{D_0(r)}{H(r)} -Cr,
\end{equation}
which shows that for every $\epsilon>0$, there exists a small $r_{\epsilon}>0$ such that $N(r)\geq -\epsilon$ for $0<r<r_{\epsilon}$. That is, \eqref{Nliminf} holds.\\

We can now proceed with the proof of Theorem \ref{Almgren's_per}.\\
\begin{pf} [Proof of Theorem \ref{Almgren's_per}]
We compute
\begin{equation}\label{Nprime}
N{'}(r)=\frac{N(r)}{r}+r\frac{D{'}(r)}{H(r)}-N(r)\frac{H{'}(r)}{H(r)}.
\end{equation}
First we see that
\begin{equation}\label{Dest}
D(r)=\int_{(\partial B_r)^+} (u u_r+ v v_r).
\end{equation}
Moreover, by direct computation, we have
\begin{equation}\label{Hprimeest}
H{'}(r)=\frac{n}{r} H(r) + 2\int_{(\partial B_r)^+} (u u_r+ v v_r).
\end{equation}
Next,
\begin{equation*}
\begin{split}
&D{'}(r)=\int_{(\partial B_r)^+}  \big(|\nabla u|^2+|\nabla v|^2 + uv\big) dS+ \int_{\partial B_r{'}} \left(\lambda_- (u^-)^{p-1}  -\lambda_+ (u^+)^{p-1}\right) vdx\\
&=\frac{n-1}{r} \int_{B_r^+}|\nabla u|^2+ 2 \int_{(\partial B_r)^+} {u_r}^2 dS -\frac{2}{r}  \int_{B_r^+} (x \cdot \nabla u) v  dx dy\\
&+ \frac{n-1}{r}\int_{B_r^+}|\nabla v|^2+ 2 \int_{(\partial B_r)^+} {v_r}^2 dS+\int_{(\partial B_r)^+} uv dS \\
&-\frac{2}{r}  \int_{B_r{'}} (x \cdot \nabla v)  \left(\lambda_- (u^-)^{p-1}  -\lambda_+ (u^+)^{p-1}\right)  dx dy+\int_{\partial B_r{'}} \left(\lambda_- (u^-)^{p-1}  -\lambda_+ (u^+)^{p-1}\right)vdx, \\
\end{split}
\end{equation*}
where we have used the Rellich identity (\ref{Rellich}) for $u$ and $v$. Therefore,
 \begin{equation}\label{estt}
\begin{split}
&D{'}(r)=\frac{n-1}{r} \int_{B_r^+}  \big(|\nabla u|^2+|\nabla v|^2 + uv\big) dx dy + 2 \int_{(\partial B_r)^+} {u_r}^2 dS+2 \int_{(\partial B_r)^+} {v_r}^2 dS\\
&+ \int_{(\partial B_r)^+} u v dS-\frac{n-1}{r} \int_{B_r^+} u v dx dy -\frac{2}{r} \int_{B_r^+}  (x\cdot \nabla u) v dx dy \\
&-\frac{2}{r}  \int_{B_r{'}} (x \cdot \nabla v)  \left(\lambda_- (u^-)^{p-1}  -\lambda_+ (u^+)^{p-1}\right)  dx dy+\int_{\partial B_r{'}} \left(\lambda_- (u^-)^{p-1}  -\lambda_+ (u^+)^{p-1}\right)vdx. \\
\end{split}
\end{equation}
Letting $$F(u):=\lambda_- (u^-)^{p-1}  -\lambda_+ (u^+)^{p-1}, $$
we have $$\nabla F(u)=G(u) \nabla u \ \  \text { where } \ \ G(u):=-(p-1)\left( \lambda_- (u^-)^{p-2} + \lambda_+ (u^+)^{p-2}\right).$$
When $p=2$, $G(u)=-\lambda_-\chi_{\{u<0\}} -\lambda_+\chi_{\{u>0\}}$ and $\nabla G(u)=0$. In the case $p\geq 3$ instead $$\nabla G(u)=-(p-1)(p-2)\left(- \lambda_- (u^-)^{p-3} + \lambda_+ (u^+)^{p-3}\right) \nabla u.$$
We now focus our attention to the last two terms in (\ref{estt}). An integration by parts on $B_r{'}$ gives
\begin{equation}\label{term}
\begin{split}
-\frac{2}{r}  \int_{B_r{'}} (x \cdot \nabla v)  F(u) dx&+\int_{\partial B_r{'}} F(u) v dx
=-\frac{1}{r}   \int_{B_r{'}} (x \cdot \nabla v)  F(u) dx-\int_{\partial B_r{'}} F(u) v dx\\&+\frac{n}{r} \int_{B_r{'}}F(u) v dx+\frac{1}{r}   \int_{B_r{'}}  G(u) v (x \cdot \nabla u)dx+ \int_{\partial B_r{'}}F(u)v dx\\
&=-\frac{1}{r}   \int_{B_r{'}} (x \cdot \nabla v)  F(u) dx+\frac{n}{r} \int_{B_r{'}}F(u) v dx+\frac{1}{r}   \int_{B_r{'}}  G(u) v (x \cdot \nabla u)dx\\
&=\frac{1}{r}   \int_{B_r{'}} x \cdot  \left( G(u) v \nabla u-F(u)\nabla v\right)  dx+\frac{n}{r} \int_{B_r{'}}F(u) v dx.
\end{split}
\end{equation}
The identity (\ref{estt}) then becomes
\begin{equation}\label{estt2}
\begin{split}
D{'}(r)&=\frac{n-1}{r}D(r) + 2 \int_{(\partial B_r)^+} {u_r}^2 dS+2 \int_{(\partial B_r)^+} {v_r}^2 dS\\
&+ \int_{(\partial B_r)^+} u v dS
-\frac{n-1}{r} \int_{B_r^+} u v dx dy -\frac{2}{r} \int_{B_r^+}  (x\cdot \nabla u) v dx dy\\&+\frac{{1}}{r} \int_{B_r{'}}F(u) v dx
+\frac{1}{r}   \int_{B_r{'}} x \cdot  \left( G(u) v \nabla u-F(u)\nabla v\right)  dx
\end{split}
\end{equation}
We need to estimate the terms in the right hand side of \eqref{estt2}. We  have
\begin{equation}\label{eq.n1}
\int_{(\partial B_r)^+} uv dS \leq C H(r)
\end{equation}
and
\begin{equation}\label{eq.n2}
\begin{split}
&-\frac{n-1}{r} \int_{B_r^+} u v dx dy \leq \frac{n-1}{r}\int_{B_r^+}\left( u^2 + v^2 \right) dx dy \leq C (H(r)+r D_0(r)),\\
\end{split}
\end{equation}
where the last inequality follows from Lemma  \ref{poincareineq}. Similarly,
\begin{equation}\label{eq.n3}
\frac{-2}{r} \int_{B_r^+}  (x\cdot \nabla u) v dx dy \leq  \int_{B_r^+} v^2 dx dy +\int_{B_r^+}  |\nabla u|^2 dx dy  \leq C(r H(r) + D_0(r) ).
\end{equation}
Next, as in \eqref{extra2} we obtain
\begin{equation}\label{eq.n4}
\frac{1}{r} \int_{B_r{'}}F(u) vdx\leq C\left(\frac{H(r)}{r}+D_0(r)\right).
\end{equation}
It remains to estimate the last term on the right hand side  of (\ref{estt2}). 
We begin by computing
\begin{equation}\label{lastterm1}
\begin{split}
\int_{B_r{'}} (x \cdot \nabla u)G(u)  v dx
&=-\int_{B_r^+} \frac{\partial}{\partial y} \left((x \cdot \nabla u) G(u) v\right) dx dy + \frac{1}{r} \int_{(\partial B_r)^+} y (x \cdot \nabla u) G(u) v dS\\
&=-\int_{B_r^+}  (x \cdot \nabla u_y) G(u) v dx dy-\int_{B_r^+}  (x \cdot \nabla u) G(u) v_y dx dy \\
& \quad -\int_{B_r^+}  (x \cdot \nabla u) \frac{d}{dy}G(u) v dx dy+ \frac{1}{r} \int_{(\partial B_r)^+}y (x \cdot \nabla u)  G(u) v dx\\
&=\int_{B_r^+}\left(n G(u) v +G(u)  x \cdot \nabla v+ v x \cdot \nabla G(u)  \right)u_y dx dy\\
&-\int_{B_r^+}  (x \cdot \nabla u) G(y) v_y dx dy
 \quad -\int_{B_r^+}  (x \cdot \nabla u) \frac{d}{dy}G(u) v dx dy\\
 &-\frac{1}{r} \int_{(\partial B_r)^+}|x|^2  v G(u)u_y dS
 \quad + \frac{1}{r} \int_{(\partial B_r)^+} y (x \cdot \nabla u)  G(u)vdx\\
 \end{split}
 \end{equation}
Taking into account Lemma \ref{poincareineq*}, we see that the first three integrals in the last equality in (\ref{lastterm1}) are bounded above by
 $C \left(r D_0(r)+H(r)\right)$.
Moreover, the last two terms in (\ref{lastterm1}) are bounded  above by $Cr \sqrt{H(r)} \sqrt{B(r)}.$ Here
\begin{equation}\label{Bdef}
B(r):=\int_{(\partial B_r)^+}\left( |\nabla u|^2 + |\nabla v|^2\right) dS.
\end{equation}
In both instances, the constant $C$ depends on the dimension and the local uniform bounds of $u$ and $v$. Similarly, we have
\begin{equation}\label{lastterm2}
\begin{split}
&\int_{B_r{'}} (x \cdot \nabla v)  F(u) dx\\
&=-\int_{B_r^+} \frac{\partial}{\partial y} \left((x \cdot \nabla_x v)  F(u)\right) dx dy + \frac{1}{r} \int_{(\partial B_r)^+} y (x \cdot \nabla_x v)  F(u) dS\\
&=-\int_{B_r^+}  (x \cdot \nabla_x v_y)  F(u) dx dy-\int_{B_r^+}  (x \cdot \nabla_x v)  G(u) u_y dx dy  + \frac{1}{r} \int_{(\partial B_r)^+}y (x \cdot \nabla_x v)  F(u) dx\\
&=\int_{B_r^+}\left(nF(u) + G(u) x \cdot \nabla_x u \right)v_y dx dy-\int_{B_r^+}  (x \cdot \nabla_x v)  G(u) u_y dx dy\\
& -\frac{1}{r} \int_{(\partial B_r)^+}|x|^2 F(u) v_y dS + \frac{1}{r} \int_{(\partial B_r)^+} y (x \cdot \nabla_x v)  F(u)dx.\\
 \end{split}
 \end{equation}
Arguing as we did for (\ref{lastterm1}), we can show that the first two integrals in  (\ref{lastterm2}) are bounded above by
 $C\left(r D_0(r)+H(r)\right)$, whereas the remaining two  are bounded  above by $Cr \sqrt{H(r)} \sqrt{B(r)}$. In both estimates, the constant $C$ depends on the dimension and the local uniform bounds of $u$ and $v$. We can thus conclude
 \begin{equation}\label{lastterm}
\frac{1}{r}   \int_{B_r{'}} x \cdot  \left( G(u) v \nabla u-F(u)\nabla v\right)  dx\leq C\left(\frac{H(r)}{r} + D_0(r)\right)+C\sqrt{H(r)}\sqrt{B(r)}.
\end{equation}
Taking into account \eqref{eq.n1}-\eqref{eq.n4} and \eqref{lastterm1}, the estimate (\ref{estt2}) becomes
 \begin{equation}\label{estt3}
 D{'}(r) \geq\frac{n-1}{r}D(r) + 2 \int_{(\partial B_r)^+} \left({u_r}^2+{v_r}^2\right) dS- C\left(\sqrt{H(r)}\sqrt{B(r)} + \frac{H(r)}{r} + D_0(r)\right),
 \end{equation}
where $B(r)$ is as defined in (\ref{Bdef}). It remains to obtain an upper bound for $B(r)$.  For this, we write
 \begin{equation*}
 \begin{split}
 D{'}(r)&= B(r) + \int_{(\partial B_r)^+} uv dS +\int_{\partial B_r{'}}F(u) v dx\\
 &=B(r) + \int_{(\partial B_r)^+} uv dS+\frac{n}{r} \int_{B_r{'}}F(u) v dx +\frac{1}{r}  \int_{B_r{'}} (x \cdot \nabla v)  F(u) dx+\frac{1}{r}   \int_{B_r{'}}  G(u) v (x \cdot \nabla u)dx\\
 &\geq B(r) -C \left(D_0(r)+\frac{H(r)}{r}\right)-C\sqrt{H(r)}\sqrt{B(r)}.\\
 \end{split}
 \end{equation*}
Here,  the first equality follows from the definition of $D(r)$, the second one  from an integration by parts on $B_r{'}$, and the last inequality from an application of \eqref{eq.n1}, \eqref{eq.n4} and \eqref{lastterm}. Hence,  by Young's inequality
\begin{equation}\label{Best1}
D{'}(r)\geq B(r)- C\frac{H(r)}{r} - CD_0(r)-\frac{1}{2} B(r),
\end{equation}
or equivalently
\begin{equation}\label{Best2}
B(r) \leq 2 D{'}(r) + C\left(\frac{H(r)}{r} +D_0(r)\right).
\end{equation}
At this point we need an upper bound for $D{'}(r)$. We consider two different cases. If $N{'}(r) \geq 0$ for a.e. $r>0$ small enough, then $N(r)$ is non-decreasing. Taking into account (\ref{Nliminf}) we conclude that $N(r)$ has a non-negative limit as  $r \to 0^+$ and the proof is complete. Thus, without loss of generality, we may now assume that for every small $r>0$, $\left|(0,r) \cap \{N{'}\leq 0\}\right|>0$. But $N{'}(r)\leq 0$ is equivalent to saying that
$$D{'}(r) \leq D(r)\left(\frac{H{'}(r)}{H(r)}-\frac{1}{r}\right)=D(r)\left(\frac{n}{r} + 2\frac{D(r)}{H(r)}-\frac{1}{r}\right)=\frac{n-1}{r} D(r) +2 \frac{D^2(r)}{H(r)}.$$
Combining this estimate with (\ref{Best2}) gives
\begin{equation}\label{Best}
B(r) \leq \frac{2(n-1)}{r} D(r) +4 \frac{D^2(r)}{H(r)} + C\left(\frac{H(r)}{r} +D_0(r)\right).
\end{equation}
Recalling \eqref{extra1} and \eqref{extra2},  for $r>0$ small enough we have
\begin{equation}\label{DtoD0}
|D(r)|\leq C\left(rD_0(r)+H(r)\right)
\end{equation}
and therefore
\begin{equation}
\begin{split}
B(r)H(r)&= \frac{2(n-1)}{r} D(r) H(r)+4 D^2(r)+ CH(r)\left(\frac{H(r)}{r} +D_0(r)\right)\\
&\leq {C} D_0(r)H(r) +CD^2_0(r) + C\frac{H^2(r)}{r}\\
&\leq CD^2_0(r) + C\frac{H^2(r)}{r},
\end{split}
\end{equation}
which yields
$$\sqrt{B(r)H(r)} \leq C\left(D_0(r)+\frac{H(r)}{r}\right).$$
Using this estimate in (\ref{estt3}) we obtain
 \begin{equation}\label{Dprimeest}
 D{'}(r) \geq\frac{n-1}{r}D(r) + 2 \int_{(\partial B_r)^+} \left({u_r}^2 +{v_r}^2\right) dS- C\frac{H(r)}{r} - C D_0(r).
 \end{equation}
Finally, using (\ref{Dest}), (\ref{Hprimeest}) and (\ref{Dprimeest}) in (\ref{Nprime}), we have
 \begin{equation*}
 \begin{split}
N{'}(r)&= \frac{N(r)}{r}+r\frac{D{'}(r)}{H(r)}-N(r)\left(\frac{n}{r} +2\frac{D(r)}{H(r)}\right)\\
&=-(n-1)\frac{N(r)}{r}+r\frac{D{'}(r) H(r) -2 D^2(r)}{H^2(r)}\\
&\geq 2r\frac{ \left( \int_{(\partial B_r)^+}
\left({u_r}^2 + {v_r}^2\right) dS\right)\left(\int_{(\partial B_r)^+} (u^2+v^2)dS\right) - \left(\int_{(\partial B_r)^+}(u u_r + v v_r)dS\right)^2}{H^2(r)}\\
&-C-Cr\frac{D_0(r)}{H(r)}.
\end{split}
\end{equation*}
Letting
\begin{equation}\label{eq.N0}
N_0(r)=r\frac{D_0(r)}{H(r)}
\end{equation}
and applying the Cauchy-Schwartz inequality we conclude
\begin{equation}\label{eq.N1}
N{'}(r)\geq -C-C N_0(r).
\end{equation}
From \eqref{N_to_N0} we infer
\begin{equation}\label{eq.equiv1}
N_0(r)\leq \frac{N(r)+Cr}{1-Cr}
\end{equation}
for some $C>0$, which combined with \eqref{eq.N1} gives
$$
N{'}(r)\geq -C-C\frac{N(r)+Cr}{1-Cr}.
$$
Thus, taking (\ref{Nliminf}) into account,  we have
\begin{equation}\label{Nmonotone}
\frac{N{'}(r)}{1+N(r)} \geq -C.
\end{equation}
We recall that we are working under the assumption $N{'}(r)<0$. But it is obvious, thanks again to (\ref{Nliminf}), that the inequality (\ref{Nmonotone}) also holds in the case $N{'}(r)\geq 0$. Hence, (\ref{Nmonotone}) holds for all $0<r<r_0$ with $r_0$ sufficiently small. For a fixed $r$ in $(0,r_0)$, we integrate (\ref{Nmonotone}) on $(r, r_0)$ to deduce
\begin{equation}\label{monotone}
e^{Cr} (N(r)+1) \leq e^{Cr_0 } (N(r_0)+1)
\end{equation}
for all $0<r<r_0$. Therefore, $e^{Cr} (N(r)+1)$ is monotone increasing. Thus, it has a limit as $r \to 0$. This along with (\ref{Nliminf}) immediately implies that $\mu:=\lim_{r \to 0} N(r)$ exists and it is non-negative.
\end{pf}

Now we observe that $N_{0}(r)$, as defined in \eqref{eq.N0}, also has a limit as $r \to 0^+$.

\begin{corollary}\label{Almgren's}
There exists $\lim_{r \to 0^+} N_{0}(r)$, and it coincides with  $\lim_{r \to 0^+} N(r)$.
\end{corollary}

\begin{pf}
Immediate consequence of the estimates (\ref{extra1}), (\ref{extra2}), and \eqref{eq.equiv1}.
\end{pf}

Next, we introduce the quantity
\begin{equation}\label{phi_def}
\phi(r) =\frac{1}{r^{n}}\int_{(\partial B_r)^+} \left(u^2+ v^2\right) dS,
\end{equation}
and we obtain the following estimates for it.

\begin{lemma} \label{phi_est}
Let $\mu:=\lim_{r \to 0^+} N(r)$. Then, the following hold:
\begin{enumerate}
\item There exists a constant $\alpha>0$, and $r_0>0$ such that
\begin{equation}\label{phi_est_upper}
r^{-2\mu} \phi(r) \leq R^{-2\mu} \phi(R)  e^{\frac{2C(R^{\alpha}-r^{\alpha})}{\alpha}}
\end{equation}
 for all $0<r<R<r_0$.
\item For any $\delta>0$, there exists $r_0(\delta)$ such that

\begin{equation}\label{phi_est_lower}
r^{-2\mu} \phi(r) \geq R^{-2\mu} \phi(R) \left(\frac{r}{R}\right)^{2\delta}
\end{equation}
for all $0<r<R<r_0(\delta)$.
\end{enumerate}
\end{lemma}
\begin{pf}
By (\ref{monotone}) we have
$N(r) \geq e^{-Cr}(\mu+1)-1$ for all $0<r<r_0$ for some $r_0>0$,
which gives
\begin{equation}\label{N_lower_bnd}
N(r)-\mu \geq (\mu+1)(e^{-Cr}-1) \geq -Cr^{\alpha }
\end{equation}
for some $\alpha>0$ and for all $0<r<r_0$. Moreover, by the definition of $\mu$,  for any $\delta>0$ there exists $r_0(\delta)$ such that
\begin{equation}\label{N_upper_bnd}
 N(r)-\mu \leq \delta
 \end{equation}
 for all $0<r<r_0(\delta)$. Now we compute
$$\phi{'}(r) =\frac{1}{r^n}\int_{(\partial B_r)^+} \left(2u u_r+ 2v v_r\right) dS=\frac{2}{r^{n}} D(r),$$
which in turn yields
\begin{equation*}
\frac{d}{dr} \left(r^{-2\mu} \phi(r)\right)=2\frac{r^{-2 \mu -1}}{r^{n}}\left(r D(r)-\mu  H(r)  \right)={2} r^{-2 \mu-1} \phi(r)\left(N(r)-\mu\right).
\end{equation*}
We can rewrite this as
\begin{equation}\label{eq.lnder}
\frac{\frac{d}{dr} (r^{-2\mu} \phi(r))}{r^{-2\mu} \phi(r)}=\frac{2}{r} (N(r)-\mu)
\end{equation}
for all $r>0$. Using the estimate (\ref{N_lower_bnd}) we get
\begin{equation*}
\frac{\frac{d}{dr} (r^{-2\mu} \phi(r))}{r^{-2\mu} \phi(r)} \geq -2C r^{\alpha-1}
\end{equation*}
for all $0<r<r_0$, and integrating over $(r,R)$ for  $0<r<R<r_0$ we obtain
\begin{equation}\label{parta}
 r^{-2\mu} \phi(r) \leq R^{-2\mu} \phi(R)  e^{\frac{2C(R^{\alpha}-r^{\alpha})}{\alpha}},
\end{equation}
which proves the first part of the lemma.

Similarly, using (\ref{N_upper_bnd}) in \eqref{eq.lnder} we get
\begin{equation*}
\frac{\frac{d}{dr} (r^{-2\mu} \phi(r))}{r^{-2\mu} \phi(r)} \leq \frac{2\delta}{r}
\end{equation*}
for all $0<r<r_0(\delta)$.  Integrating over $(r,R)$ for  $0<r<R<r_0(\delta)$ we now infer
\begin{equation}\label{partb}
r^{-2\mu} \phi(r) \geq R^{-2\mu} \phi(R) \left(\frac{r}{R}\right)^{2\delta},
\end{equation}
thus showing the second part of the lemma.
\end{pf}

As a corollary, we obtain the following growth estimate.
\begin{corollary} \label{uv_avg}
Let $u$ be the minimizer of (\ref{fnct}). Then, for all $z \in B^+_r\cup {B_r^{'}}$, $0<r<\frac{1}{2}$, we have
$$u(z) \leq C r^{\mu} \ \  \text { and } \ \ v(z) \leq C r^{\mu}, $$
where $C$ depends on the dimension and the local $\mathcal{L}^{\infty}$-norms of $u$ and $v$.
\end{corollary}
\begin{pf}
Since $u^{\pm}$ are subharmonic, we have
\begin{equation}\label{u_growth}
u^{\pm}(z) \leq \intbar_{(\partial B_{r})^+} u^{\pm} dS<Cr^{\mu},
\end{equation}
where the last inequality follows from (\ref{phi_est_upper}). The constant $C$ depends on the dimension and $||u||_{\mathcal{L}^{\infty}(B_{\frac{1}{2}})}$. Since $v^{\pm}$ are also subharmonic, we see that the same growth estimate holds for $v$ as well.
\end{pf}

 \section{Almgren rescalings and blow-ups}\label{blow_up_analysis}
 First we note that all the results of Section \ref{Almgren's_monotonicity} obtained for the origin, apply to any free boundary point $z_0=(x_0,0)$ as well.

The next step in our program is the study of blow-up sequences around a free boundary point $z_0=(x_0,0)$. Without loss of generality, we assume $z_0=0$ and let $\mu:=\lim_{r \to 0^+} N(r)$. For $z=(x,y) \in B_1^+\cup B_1^{'}$ and $0<r<1$, we define the \emph{Almgren's rescalings}
\begin{equation}\label{Almgren's_rescaling}
u_r(z):= \frac{u(rz)}{\sqrt{\phi(r)}} \text { and } v_r(z):= \frac{ v(rz)}{\sqrt{\phi(r)}}
\end{equation}
We notice that $\Delta u_r= r^2 v_r$. In addition,  we have
\begin{equation}\label{normalization}
||u_r||_{\mathcal{L}^2(\partial B_1)^+}+ ||v_r||_{\mathcal{L}^2(\partial B_1)^+} =1
\end{equation}
The following convergence theorem is crucial in our analysis.
\begin{theorem}\label{blow-up}
Let $u$ be the solution to (\ref{prb}) with $p=2$ or $p\geq 3$, and $v=\Delta u$. There exists subsequences $u_j:=u_{r_j}$ and $v_j:=v_{r_j}$ such that $u_j \to u^{\ast}$, $\nabla u_j \to \nabla u^{\ast}$, $v_j \to v^{\ast}$ and $\nabla v_j \to \nabla v^{\ast}$ as $j\to\infty$, uniformly on every compact subset of $\overline{\mathbb{R}_+^{n+1}}$. Moreover, $u^{\ast}$ and $v^{\ast}$ are homogenous harmonic polynomials of degree $\mu$. If, in addition, we have $\nabla_x u(0)=0$ or $\nabla_x v(0)=0$, then we must have $\mu \geq 2$.
\end{theorem}

\begin{pf}
First we want to show that the sequence $\{u_r\}$ is bounded in $W_{loc}^{2,2}(\overline{\mathbb{R}_+^{n+1}})$. To this end, we begin by noticing that for $r, R>0$ we have
\begin{equation}\label{N_r}
\begin{split}
N(rR,u,v)=&rR\frac{\int_{B_{rR}^+} \left(|\nabla u|^2+|\nabla (\Delta u)|^2 + u\Delta u\right) dx dy + \int_{B_{rR}{'}} \left(\lambda_- (u^-)^{p-1}  -\lambda_+ (u^+)^{p-1}\right)v dx}{\int_{(\partial B_{rR})^+} (u^2+ (\Delta u)^2) dS}\\
&\geq R \frac{\int_{B_{R}^+} \left(|\nabla u_r|^2+r^{-4}|\nabla (\Delta u_r)|^2 + u_r\Delta u_r\right) dx dy -C r^{-1} \int_{B_R{'}} |u_r|
|\Delta u_r|dx}{\int_{(\partial B_R)^+} (u_r^2+ r^{-4} (\Delta u_r)^2) dS}\\
&=R\frac{\int_{B_{R}^+} \left(|\nabla u_r|^2+|\nabla v_r|^2 +r^2 u_r v_r\right) dx dy -Cr \int_{B_R{'}} |u_r|
|v_r|dx}{\int_{(\partial B_R)^+} (u_r^2+  v_r^2) dS}.\\
\end{split}
\end{equation}
Here, the constant $C$ depends on $\lambda_{\pm}$, $p$, and the local $\mathcal{L}^{\infty}$-bound of $u$.
Thus,
\begin{equation*}
\begin{split}
&\int_{B_{R}^+} (|\nabla u_r|^2+|\nabla v_r|^2) dx dy\\
 &=\frac{1}{R}N(rR,u,v)\left(\int_{(\partial B_R)^+}\left(u_r^2 + v_r^2 \right)dS\right)-r^2 \int_{B_{R}^+}u_r v_r dx dy + C r \int_{B_R{'}}\left(u_r^2 + v_r^2 \right)dx\\
&\leq\frac{1}{R} N(rR,u,v)\left(\int_{(\partial B_R)^+}\left(u_r^2 + v_r^2 \right)dS\right) + \frac{r^2}{2}  \int_{B_{R}^+} \left( u_r^2 + v_r^2\right)dx dy+ C r \int_{B_R{'}}\left(u_r^2 + v_r^2 \right)dx\\
&\leq\frac{1}{R} N(rR,u,v) \left(\int_{(\partial B_R)^+}\left(u_r^2 + v_r^2 \right)dS\right)+ C r^2 \int_{B_{R}^+} (|\nabla u_r|^2+|\nabla v_r|^2) dx dy+Cr\int_{(\partial B_R)^+}\left(u_r^2 + v_r^2 \right)dS\\
\end{split}
\end{equation*}
where, in the last inequality, we have used Lemma \ref{poincareineq*} and Lemma \ref{l.trace}. Therefore, for $R_0(\delta)$ as in Lemma \ref{phi_est} and $r>0$ small enough such that $rR<R_0(\delta)$  and $Cr^2<1$, we have
\begin{equation}\label{almg_uni_bnd}
\begin{split}
\int_{B_{R}^+} (|\nabla u_r|^2+|\nabla v_r|^2) dx dy& \leq C  \left(\frac{1}{R}N(rR,u,v)+r\right)\int_{(\partial B_R)^+}\left(u_r^2 + v_r^2 \right)dS\\
&=C \left(\frac{1}{R}N(rR,u,v)+r\right)\frac{\phi(rR,u)}{\phi(r,u)}R^n\\
&\leq C \left((\mu+\delta)+rR\right) R^{n-1+2(\mu+\delta)}.
\end{split}
\end{equation}
The last inequality follows from (\ref{phi_est_lower}) and (\ref{N_upper_bnd}). We conclude that, for every fixed $R>1$, the sequences $\{u_r\}$ and $\{v_r\}$ are uniformly bounded in $W^{1,2}_{loc}(B_R^+\cup B_R^{'})$ for $r$ sufficiently small. Moreover, it is readily seen that the pair $(u_r, v_r)$ is a weak solution of the system
\begin{equation}\label{rescalingsys}
\left\{
\begin{array}{ll}
\Delta u_r =r^2 v_r & \text { in } B_R^+\\
\Delta v_r =0 & \text { in } B_R^+\\
(u_r) _y=0 & \text { on } B_R{'}\\
(v_r)_{y}=r \left(\phi(r,u)\right)^{\frac{p-2}{2}}\left(\lambda_- (u_r^-)^{p-1}-\lambda_+ (u_r^+)^{p-1} \right)& \text { on }  B_R{'}\\
\end{array}
\right.
\end{equation}

In a similar way to the proof of Theorem \ref{part_reg}, we can show that $u_r \in C_{loc}^{3,\alpha}(\overline{\mathbb{R}_+^{n+1}})$   and $v_r \in C_{loc}^{1,\alpha}(\overline{\mathbb{R}_+^{n+1}})$ with norms controlled by $||u_r||_{W^{2,2}(B_R^+)}+||v_r||_{W^{1,2}(B^+_R)}$. This result, along with a compactness argument, implies that for any sequence $\{(u_{r_j}, v_{r_j})\}$, with $r_j \to 0^+$ as $j \to \infty$, there exists a subsequence denoted by $\{(u_r, v_r)\}$ and two functions $ u^{\ast}$ and $ v^{\ast}$ such that
$u_j \to u^{\ast}$, $\nabla u_j \to \nabla u^{\ast}$, $\Delta u_j \to \Delta u^{\ast}$,  $v_j \to v^{\ast}$ and $\nabla v_j \to \nabla v^{\ast}$ uniformly on every compact subset of $\mathbb{R}^{n+1}$.  Moreover, it is obvious  that
\begin{equation}
||u^{\ast}||_{\mathcal{L}^2(\partial B_1^+)}+ ||v^{\ast}||_{\mathcal{L}^2(\partial B_1^+)} =1,
\end{equation}
which implies that $u^{\ast}$ and  $v^{\ast}$ cannot both vanish identically. Next, passing to the weak limit in (\ref{rescalingsys}) with $r=r_j \to 0^+$, we see that  $u^{\ast}$ and $v^{\ast}$ are both harmonic in $B_R^+$ and satisfy $u^{\ast}_y=v^{\ast}_y=0$ on $B_R{'}$. Furthermore, we have
\begin{equation*}
N_0(R , u^{\ast}, v^{\ast})=\lim_{j \to 0} N_0( R, u_{r_j}, v_{r_j})=\lim_{j \to \infty} N_0(r_j R,u,v)=\lim_{j \to \infty} N(r_j R,u,v)=\mu.
\end{equation*}
Here, the second equality is obtained by rescaling, whereas the last two follow from Corollary \ref{Almgren's} and Theorem \ref{Almgren's_per}, respectively. Arguing as in the proof of \cite[Lemma 3.1]{FF}, one can see that  $\mu$ must be an integer, and $u^{\ast}$ and $v^{\ast}$ must be homogenous harmonic polynomials of degree $\mu$.

Finally, if in addition we have $\nabla_x u(0)=0$ or $\nabla_x v(0)=0$, then, by the uniform convergence of $\nabla u_j$ and $\nabla v_j$ to $\nabla u^{\ast} $ and $\nabla v^{\ast}$ respectively, we must have $\mu \geq 2$.
\end{pf}

\section{A Monneau-type monotonicity formula}\label{Monneau_monotonicity}
Now we turn our attention to the study of the regularity and structure of the free boundary. We define the \emph{regular set} of the free boundary as
\begin{equation}\label{regularset}
\mathcal{R}(u)=\{z=(x,0) \in \Gamma(u) \ \ | \ \  u(x,0)=0, \nabla_x u (x,0) \neq 0, \text { and } \nabla_x v (x,0) \neq 0 \}
\end{equation}
and the\emph{ singular set}  $\Sigma(u):=\Gamma (u) \setminus \mathcal{R}(u)$. An immediate consequence of the regularity of the solution and of the implicit function theorem is the following  result showing the regularity of the free boundary around  regular points.
\begin{theorem}\label{regular_fb}
If $(x,0) \in \mathcal{R}(u)$, then the free boundary is a $C^{3,\alpha}-$ graph in a neighborhood of $(x,0)$  for all $0<\alpha<1$.
\end{theorem}

To study the structure of the singular set, we first need to derive some additional properties (such as non-degeneracy and  continuous dependence of blow-ups) at singular points. A crucial tool in this endeavor is the monotonicity of the Monneau-type functional
 \begin{equation}\label{monneau}
M_{\mu} (r,u, v, p_{\mu}, q_{\mu})=\frac{1}{r^{n+2\mu}} \int_{(\partial B_r)^+} \left((u-p_{\mu})^2 + (v-q_{\mu})^2\right) dS
\end{equation}
 where $p_{\mu}$ and $q_{\mu}$ are two homogeneous harmonic polynomials of degree $\mu$. We introduce the Weiss-type functional
\begin{equation}\label{weiss}
W_{\mu} (r,u,v)=\frac{H(r,u,v)}{r^{n+2\mu}} (N_0(r,u,v)-\mu)
\end{equation}
to prove the following result.

\begin{theorem}\label{thm_monneau}
Let $u$ be the solution to (\ref{prb}), and let $p_{\mu}$ and $q_{\mu}$ be two homogenous harmonic polynomials of degree $\mu$, symmetric with respect to the thin space $\{y=0\}$. Then, there exists a positive constant $C$ such that

\begin{equation}\label{monneauest1}
\frac{d}{dr} M_{\mu} (r,u, v, p_{\mu}, q_{\mu}) \geq \frac{2}{r}W_{\mu} (r,u,v)-C
\end{equation}
\end{theorem}

\begin{pf}
We begin by observing that, for all $0<r<1$, $N_0(r, p_{\mu} , q_{\mu} )=\mu$ and thus  $W_{\mu}(r,p_{\mu},q_{\mu})=0$. We let $w=u-p_{\mu}$ and $h=v-q_{\mu}$. Then,
\begin{equation*}
\begin{split}
W_{\mu} (r,u,v)&=W_{\mu} (r,w+p_{\mu}, h+q_{\mu})\\
&=\frac{1}{r^{n+2\mu-1}} D_0(r,w+p_{\mu}, h+ q_{\mu})-\frac{\mu}{r^{n+2\mu}} H(r,w+p_{\mu}, h+ q_{\mu})\\
&=\frac{1}{r^{n+2\mu-1}} \int_{B_r^+} \left(|\nabla w|^2+|\nabla h|^2 + 2 \nabla w \cdot \nabla p_{\mu}+ 2 \nabla h \cdot \nabla q_{\mu}\right) dx dy\\
&- \frac{\mu}{r^{n+2\mu}} \int_{(\partial B_r)^+} \left(w^2 + h^2+ 2 w p_{\mu}+ 2 h q_{\mu}\right) dS.\\
\end{split}
\end{equation*}
Integrating by parts and keeping in mind that $p_{\mu}$ and $q_{\mu}$ are both harmonic, with $(p_{\mu})_y=(q_{\mu})_y=0$ on $B_r{'}$, we have
\begin{equation}\label{w_est1}
\begin{split}
W_{\mu} (r,u,v)&=\frac{1}{r^{n+2\mu-1}} \int_{B_r^+} \left(-w \Delta u-h \Delta v\right) dx dy\\
&+\frac{1}{r^{n+2\mu-1}} \int_{(\partial B_r)^+} \left(w \nabla w\cdot \frac{z}{r}+h \nabla h \cdot \frac{z}{r} + 2 w  \nabla p_{\mu}\cdot \frac{z}{r} + 2 h  \nabla q_{\mu} \cdot \frac{z}{r}\right) dS\\
&-\frac{\mu}{r^{n+2\mu}} \int_{(\partial B_r)^+} \left(w^2 + h^2+ 2 w p_{\mu}+ 2 h q_{\mu}\right) dS-\frac{1}{r^{n+2\mu-1}} \int_{B_r{'}}\left( w u_y+h v_y\right) dx\\
&=-\frac{1}{r^{n+2\mu-1}} \int_{B_r^+} w v dx dy+\frac{1}{r^{n+2\mu-1}} \int_{(\partial B_r)^+} \left(w \nabla w+h \nabla h\right) \cdot \frac{z}{r} dS\\
&-\frac{\mu}{r^{n+2\mu}} \int_{(\partial B_r)^+} \left(w^2 + h^2\right) dS-\frac{1}{r^{n+2\mu-1}} \int_{B_r{'}} h v_y dx.\\
\end{split}
\end{equation}
Here  we have used $\nabla p_{\mu} \cdot z= \mu p_{\mu}$ and $\nabla q_{\mu} \cdot z= \mu q_{\mu}$. Now we notice
\begin{equation}\label{terms23}
\begin{split}
\frac{1}{r^{n+2\mu-1}} \int_{(\partial B_r)^+} \left(w \nabla w+h \nabla h\right) \cdot \frac{z}{r} dS&-\frac{\mu}{r^{n+2\mu}} \int_{(\partial B_r)^+} \left( w^2 + h^2\right) dS\\
&=\frac{r}{2}\frac{d}{dr}\left( \frac{1}{r^{n+2\mu}}\int_{(\partial B_r)^+} \left(w^2 + h^2\right) dS \right)=\frac{r}{2} \frac{d}{dr} M_{\mu} (r,u, v, p_{\mu}, q_{\mu}).
\end{split}
\end{equation}
Moreover,
\begin{equation}\label{terms14}
\begin{split}
&-\frac{1}{r^{n+2\mu-1}} \int_{B_r^+}wv\ dx dy-\frac{1}{r^{n+2\mu-1}} \int_{B_r{'}} h v_y\ dx\\
&=-\frac{1}{r^{n+2\mu-1}} \int_{B_r^+}v(u-p_{\mu})dx dy-\frac{1}{r^{n+2\mu-1}} \int_{B_r{'}}\left( \lambda_- (u^-)^{p-1}-\lambda_+ (u^+)^{p-1} \right)  (v-q_{\mu})dx\\
&\leq \frac{C}{r^{n+2\mu-1}} r^{2 \mu+n}\leq Cr,\\
\end{split}
\end{equation}
where in the first inequality above we have used Corollary \ref{uv_avg} and the fact that $p_{\mu}$ and $q_{\mu}$ are homogenous polynomials of degree $\mu$. The constant $C$ depends on the coefficients of the polynomials $p_{\mu}$ and $q_{\mu}$, and on the local $\mathcal{L}^{\infty}$-norms of $u$ and $v$. Inserting \eqref{terms23} and \eqref{terms14} in \eqref{w_est1} we finally obtain
$$W_{\mu} (r,u,v) \leq \frac{r}{2}\frac{d}{d r} M_{\mu}(r,u, v, p_{\mu}, q_{\mu}) + Cr,$$
which yields the desired conclusion.
\end{pf}

In what follows, we will always let $\mu:=\lim_{r \to 0^+} N(r)$.

\begin{corollary}\label{cor_Monneau}
Let $u$, $v$, $p_{\mu}$ and $q_{\mu}$ be as in Theorem \ref{thm_monneau}. Then, there exists a positive constant $C$ such that
$$\frac{d}{dr} M_{\mu}(r,u, v, p_{\mu}, q_{\mu}) \geq -C.$$ In particular, $M_{\mu}(r,u, v, p_{\mu}, q_{\mu})$ has a limit as $r \to 0^+$.
\end{corollary}
\begin{pf}
From (\ref{extra1}) and (\ref{extra2}) we know that
\begin{equation}\label{est1}
N_0(r,u,v) \geq \frac{N(r,u,v)-Cr}{1+Cr},
\end{equation}
and thus, applying (\ref{N_lower_bnd}),
\begin{equation*}
N_0(r,u,v)-\mu \geq \frac{-Cr}{1+Cr}.
\end{equation*}
But then, thanks to (\ref{phi_est_upper}),
$$W_{\mu} (r,u,v) \geq -C\frac{H(r,u,v)}{r^{n+2\mu}(1+Cr)} r\geq -Cr. $$
Using this in (\ref{monneauest1}) we obtain
\begin{equation*}
\frac{d}{dr} M_{\mu} (r,u, v, p_{\mu}, q_{\mu}) \geq -C,
\end{equation*}
and the proof is complete.
\end{pf}

We now prove the non-degeneracy of the solution at singular free boundary points.

\begin{lemma}\label{non_degen}
Let $u$ be the solution to (\ref{prb}) with $u(0)=0$ and either $\nabla _{x} u(0)=0$ or $\nabla _{x} v(0)=0$. There exists $C>0$ and $0<R_0<1$, depending possibly on $u$, such that
\begin{equation}\label{non_deg}
\sup_{(\partial B_r)^+} |u| \geq C r^{\mu}\text { or } \sup_{(\partial B_r)^+} |v| \geq C r^{\mu}
\end{equation}
for all $0<r<R_0$.
\end{lemma}
\begin{pf}
We argue by contradiction. Suppose (\ref{non_deg}) does not hold. Then, there exists a sequence of radii $r_j \searrow 0$  and a sequence of positive constants $C_j \searrow 0$ such that
\begin{equation}\label{assumption}
 \frac{1}{r_j^{\mu}}\sup_{(\partial B_{r_j})^+} |u|\leq C_j \text { and } \frac{1}{r_j^{\mu}}\sup_{(\partial B_{r_j})^+} |v|\leq C_j.
  \end{equation}
This implies that both
$\frac{1}{r_j^{n+2\mu}} \int_{(\partial B_{r_j})^+}u^2$ and $\frac{1}{r_j^{n+2\mu}} \int_{(\partial B_{r_j})^+}v^2$ converge to $0$ as $r_j \to 0^+$. We know by Theorem \ref{blow-up} that the Almgren rescalings  $u_r$ and $v_r$ (as defined in \eqref{Almgren's_rescaling}) converge - passing to a subsequence $r_j$ - to some homogenous harmonic polynomials $p_\mu$ and $q_\mu$ of degree $\mu\geq 2$. We now  compute
\begin{equation}\label{M_est}
\begin{split}
M_{\mu} (r, u, v, p_{\mu}, q_{\mu})&=\frac{1}{r^{n+2\mu}} \int_{(\partial B_r)^+}(u^2 + v^2) dS\\
&+\frac{1}{r^{n+2\mu}} \int_{(\partial B_r)^+}(-2 up_{\mu}-2vq_{\mu})dS+\frac{1}{r^{n+2\mu}} \int_{(\partial B_r)^+}(p_{\mu}^2 + q_{\mu}^2) dS.
\end{split}
\end{equation}
By the assumption (\ref{assumption}), the first term in the right-hand-side of (\ref{M_est}) converges to $0$ as $r=r_j \to 0^+$. For the third term on the right-hand-side of (\ref{M_est}), we notice that by the homogeneity of $p_{\mu}$ and $q_{\mu}$ we have
\begin{equation}\label{homo_poly}
\frac{1}{r^{n+2\mu}} \int_{(\partial B_r)^+}(p_{\mu}^2 + q_{\mu}^2) dS = \int_{(\partial B_1)^+}(p_{\mu}^2 + q_{\mu}^2) dS.
\end{equation}
Finally, the second term in the right-hand-side of (\ref{M_est}) can be estimated as
\begin{equation*}
\begin{split}
&\left|\frac{1}{r^{n+2\mu}} \int_{(\partial B_r)^+}(-2 up_{\mu}-2vq_{\mu}) dS\right|\\
& \leq 2\left(\frac{1}{r^{n+2\mu}} \int_{(\partial B_r)^+}u^2\right)^{\frac{1}{2}}\left(\frac{1}{r^{n+2\mu}} \int_{(\partial B_r)^+}p_{\mu}^2\right)^{\frac{1}{2}}+ 2\left(\frac{1}{r^{n+2\mu}} \int_{(\partial B_r)^+}v^2\right)^{\frac{1}{2}}\left(\frac{1}{r^{n+2\mu}} \int_{(\partial B_r)^+}q_{\mu}^2\right)^{\frac{1}{2}}\\
\end{split}
\end{equation*}
which, again  by assumption (\ref{assumption}), converges to zero as $r=r_j \to 0^+$. Therefore, passing to the limit in (\ref{M_est}) as $r=r_j \to 0^+$ we infer
\begin{equation}\label{M_limt}
M(0+, u, v, p_{\mu}, q_{\mu})=\int_{(\partial B_1)^+}(p_{\mu}^2 + q_{\mu}^2) dS.
\end{equation}
Applying Corollary \ref{cor_Monneau} we obtain
\begin{equation*}
M_{\mu} (0+, u, v, p_{\mu}, q_{\mu}) \leq M_{\mu} (r, u, v, p_{\mu}, q_{\mu}) + Cr
\end{equation*}
for all $0<r<\frac{1}{2}$, which immediately gives, together with (\ref{homo_poly}),
\begin{equation*}
\frac{1}{r^{2\mu}} \int_{(\partial B_1)^+}\left(u^2(rz) + v^2(rz)\right) dS-2\frac{1}{r^{\mu}} \int_{(\partial B_1)^+} \left(u (rz) p_{\mu}(x)+v(rz) q_{\mu}(z)\right) dS \geq -Cr
\end{equation*}
for all $0<r<\frac{1}{2}$.  We can rewrite this as
\begin{equation*}
\frac{1}{r^{2\mu}} \int_{(\partial B_1)^+}\left(\phi(r) u_r^2(z) +\phi(r)  v_r^2(z)\right) dS-2\frac{1}{r^{\mu}} \int_{(\partial B_1)^+} \left(\sqrt{\phi(r)} u _r(z) p_{\mu}(z)+\sqrt{\phi(r)} v_r(z) q_{\mu}(z)\right) dS \geq -Cr,
\end{equation*}
which is equivalent to
\begin{equation*}
\begin{split}
 \int_{(\partial B_1)^+}\left(\frac{\sqrt{\phi(r)}}{r^{\mu}} u_r^2(z) +\frac{\sqrt{\phi(r)}}{r^{\mu}} v_r^2(z)\right) dS-2 \int_{(\partial B_1)^+} \left( u _r(z) p_{\mu}(x)+v_r(z) q_{\mu}(z)\right)dS &\geq -C\frac{r^{\mu+1}}{\sqrt{\phi(r)}}\\
 &\geq -Cr^{1-\delta}
 \end{split}
\end{equation*}
for all $0<r<\frac{1}{2}$. Here,  we have employed Lemma \ref{phi_est}, part 2. Finally, letting $r=r_j \to 0^+$ in the above inequality we obtain
\begin{equation*}
- \int_{(\partial B_1)^+}  \left(p_{\mu}^2(x)+q^2_{\mu}(x)\right) dS\geq 0,
\end{equation*}
which gives a contradiction since $p_{\mu}$ and $q_{\mu}$ do not both vanish at the same time.
\end{pf}

\section{The structure of the singular set}
\label{structure_singular_set}
Let $z_0$ be a point on the free boundary of $u$, the solution of (\ref{prb}).
For $R\geq 1$, we introduce the \emph{homogenous rescalings}
\begin{equation}\label{homo_rescaling}
u_r^{\mu}(z):=\frac{u(z_0+rz)}{r^{\mu}} \text { and } v_r^{\mu}(z):=\frac{v(z_0+rz)}{r^{\mu}}, \text { for  } z \in B_R^+ \cup B^{'}_R\text { and } 0<r<\frac{1}{R},
\end{equation}
and we study the existence and uniqueness of their blow-ups.
\begin{theorem}\label{homo_blow-up}
Let $u$ be the solution of (\ref{prb}) and suppose that $z_0$ is a free boundary point such that $\nabla _{x}u(z_0)~=~0$ or $\nabla _{x}v(z_0)=0$. Then, there exist two unique harmonic polynomials $p_\mu$ and $q_\mu$, homogenous of degree $\mu$  and such that $u_r^{\mu} \to p_\mu$  in $W_{loc}^{3,2}(B^+_R \cup B_R^{'})$ and in $C^{3, \alpha}(B_R^+ \cup B_R^{'})$, and $v_r^{\mu} \to q_\mu$ in $W_{loc}^{1,2}(B^+_R)$ and in $C^{1, \alpha}(B_R^+\cup B_R^{'})$ for all $R\geq 1$.
\end{theorem}

\begin{pf}
Without loss of generality, we assume $z_0=0$.
By (\ref{phi_est_upper}) and (\ref{non_deg}) there exist positive constants $C_1$ and $C_2$  such that
\begin{equation*}
C_1 r^{\mu} \leq (\phi(r))^{\frac{1}{2}} \leq C_2 r^{\mu} \text { for all } 0<r<\frac{1}{2}.
\end{equation*}
Hence
\begin{equation*}
C_2 u_r^{\mu}(z) \leq u_r(z) \leq C_1 u_r^{\mu}(z) \text { and } C_2 v_r^{\mu}(z) \leq v_r(z) \leq C_1 v_r^{\mu}(z)
\end{equation*}
where $u_r$ and $v_r$ are the Almgren's rescalings defined in (\ref{Almgren's_rescaling}). This, together with  Theorem \ref{blow-up}, implies that for any sequence $r_j \to 0^+$ there exists a subsequence, which we will not rename for simplicity, such that $u^{\mu}_{r_j}$ and $v^{\mu}_{r_j}$ converge. Moreover, the limit functions are homogenous harmonic polynomials of degree $\mu$. It remains to be shown that all subsequences converge to the same limit. For, suppose that $(u_{\ast}, v_{\ast})$ and $(u_0, v_0)$ are limits over two different subsequences. By Corollary \ref{cor_Monneau},

\begin{equation}\label{M_limit}
\begin{split}
M_{\mu} (0+, u, v, u_{\ast}, v_{\ast})=\lim_{r_j \to 0} M_{\mu} (r_j, u, v, u_{\ast}, v_{\ast}) &=\lim_{r_j \to 0} M_{\mu} (1, u^{\mu}_{r_j}, v^{\mu}_{r_j}, u_{\ast}, v_{\ast})\\
&= \lim_{r_j \to 0}\int_{(\partial B_1)^+} (u_{r_j}-u_{\ast})^2 + (v_{r_j}-v_{\ast})^2 dS=0.
\end{split}
\end{equation}
In a similar fashion we have  $$ M_{\mu} (0+, u, v, u_0, v_0)=0$$
and therefore $$\int_{(\partial B_1)^+} (u_{\ast}-u_0)^2 + (v_{\ast}-v_0)^2 dS=0,$$
which implies that $u_0=u_{\ast}$ and $v_0=v_{\ast}$.
\end{pf}

In the remainder of the paper,\ $N_0^{z_0}$ will denote Almgren's frequency formula as in (\ref{Almgren's_frq}) at the free boundary point $z_0$,\  $\mu_{z_0}=\lim_{r \to 0} N_0^{z_0}$, and
$p_{\mu}^{z_0}$ and $q_{\mu}^{z_0}$ the $\mu-$homogenous blow-up polynomials at $z_0$ as in Theorem \ref{homo_blow-up}.
For $\mu \in \mathbb{N}$,  we define
\begin{equation}\label{mu_singular}
\Sigma_{\mu}(u)=\{(x,0) \in \Gamma(u) \ \  | \ \  u(x,0)=0,  \nabla_x u (x,0)=0 \text { \emph{or} } \nabla_x u (x,0)=0,  \text { and } \ \  N_0^{z_0} (0+,u,v)=\mu\}.
\end{equation}
We denote by $\mathcal{P}(u)$ the class of homogenous harmonic polynomials of degree $\mu$ and even in $y$. We note that $\mathcal{P}(u)$ is a convex subset of the finite dimensional vector space of all polynomials homogenous of degree $\mu$. We endow  $\mathcal{P}(u)$ with the norm of  $\mathcal{L}^2((\partial B_1)^+)$. Our next step consists in establishing  continuous dependence of the blow-ups on the singular points.

\begin{theorem}\label{cont_dep}
Let $u$ be the solution to (\ref{prb}). The mapping $z_0 \mapsto (p_{\mu}^{z_0}, q_{\mu}^{z_0})$ from $\Sigma_{\mu}(u)$ to $\mathcal{P}(u) \times \mathcal{P}(u)$  is continuous. Moreover, for any compact set $K \subset \subset \Sigma_{\mu}(u)$, there exists a modulus of continuity $\omega_\mu$ such that
$$|u(z)-p^{z_0}_{\mu}(z-z_0)|\leq \omega_\mu (|z-z_0|) |z-z_0|^{\mu}$$
for any $z_0 \in K$.
\end{theorem}

\begin{pf}
Fix $z_0 \in \Sigma_{\mu}(u)$ and $\epsilon>0$ small enough. Then, by (\ref{M_limit}),  there exists $r_{\epsilon}=r_{\epsilon}(z_0)>0$ such that
 \begin{equation}\label{x0est}
M^{z_0}_{\mu} (r_{\epsilon},u, v, p^{z_0}_{\mu}, q^{z_0}_{\mu}):=\frac{1}{r_{\epsilon}^{n+2\mu}} \int_{(\partial B_{r_\epsilon})^+} \left(\left(u(z_0+z)-p^{z_0}_{\mu}\right)^2 + \left(v(z_0+z)-q^{z_0}_{\mu}\right)^2\right) dS<\epsilon.
\end{equation}
Also, by the continuity of $u$ and $v$, there exists $\rho_{\epsilon}>0$ such that for all $z_1 \in \Sigma_{\mu}(u) \cap B_{\rho_{\epsilon}}(z_0)$  we have
 \begin{equation*}
M^{z_1}_{\mu} (r_{\epsilon},u, v, p^{z_0}_{\mu}, q^{z_0}_{\mu})<2\epsilon.
\end{equation*}
Moreover, by Corollary \ref{cor_Monneau} and possibly for a smaller $r_{\epsilon}$ we have
\begin{equation}\label{x1est}
M^{z_1}_{\mu} (r,u, v, p^{z_0}_{\mu}, q^{z_0}_{\mu})<3\epsilon
\end{equation}
for all $0<r<r_{\epsilon}$. Rescaling and passing to the limit as $r\to 0$ we see that
\begin{equation}\label{map_cont}
\int_{(\partial B_1)^+} \left(p^{z_1}_{\mu}-p^{z_0}_{\mu}\right)^2 + \left(q^{z_1}_{\mu}-q^{z_0}_{\mu}\right)^2=M^{z_1}_{\mu} (0+,u, v, p^{z_0}_{\mu}, q^{z_0}_{\mu})<3\epsilon,
\end{equation}
and this proves the continuity of the mapping.

For the second part of the theorem, we notice  that (\ref{x1est}) and (\ref{map_cont}) imply that
\begin{equation*}
M^{z_1}_{\mu} (r,u, v, p^{z_1}_{\mu}, q^{z_1}_{\mu}) \leq M^{z_1}_{\mu} (r,u, v, p^{z_0}_{\mu}, q^{z_0}_{\mu})+\frac{1}{r^{n+2\mu}}\int_{(\partial B_r)^+} \left(p^{z_1}_{\mu}-p^{z_0}_{\mu}\right)^2 + \left(q^{z_1}_{\mu}-q^{z_0}_{\mu}\right)^2 \leq 6 \epsilon
\end{equation*}
for all $x_1 \in B_{\rho_{\epsilon}} (x_0)$ and all $0<r<r_{\epsilon}$.
This is equivalent to say
\begin{equation}\label{L2_uni_bnd}
||u^{\mu}_{z_1,r}-p_{\mu}||^2_{\mathcal{L}^2(\partial B_1)^+}+||v^{\mu}_{z_1,r}-q_{\mu}||^2_{\mathcal{L}^2(\partial B_1)^+}  \leq 6\epsilon.
\end{equation}
Now we observe that $w^{\mu}_r:=u^{\mu}_r-p_{\mu}$ satisfies the oblique derivative problem
\begin{equation*}
\left\{
\begin{array}{ll}
\Delta w^{\mu}_r =r^2 v^{\mu}_r & \text { in } B_1^+\\
(w^{\mu}_r) _y=0 & \text { on } B_1'\\
\end{array}
\right.
\end{equation*}
By the  $L^{\infty}-L^2$ interior estimate for oblique derivatives problems for elliptic equations (see for example \cite{L}, theorem 5.36), we have
\begin{equation}\label{Linfty_uni_bnd_1}
||u^{\mu}_{z_1,r}-p_{\mu}||_{\mathcal{L}^{\infty} (B^+_{1/2})}
\leq \left( ||u^{\mu}_{z_1,r}-p_{\mu}||_{\mathcal{L}^2 (B^+_1)}+C r^2 ||v_r^{\mu}||_{\mathcal{L}^{\infty}(B^+_1)}\right) \leq C(\epsilon)
\end{equation}
for some $C(\epsilon) \to 0$ as $\epsilon \to 0$. Here,  the last inequality follows from (\ref{L2_uni_bnd}) and Corollary \ref{uv_avg}. Similarly, $h^{\mu}_r:=v^{\mu}_r-q_{\mu}$
is a weak solution for the oblique problem
 \begin{equation*}
\left\{
\begin{array}{ll}
\Delta h^{\mu}_r =0 & \text { in } B_1^+\\
(h^{\mu}_r)_{y}=r \left(\lambda_-\left( (u_r^{\mu})^{-}\right)^{p-1}-\lambda_+ \left((u_r^{\mu})^{+}\right)^{p-1} \right)& \text { on }  B_1'\\
\end{array}
\right.
\end{equation*}
and the same $L^{\infty}-L^2$ interior estimate implies that
\begin{equation}\label{Linfty_uni_bnd_2}
||v^{\mu}_{z_1,r}-q_{\mu}||_{\mathcal{L}^{\infty}(B^+_{1/2})} \leq C\left(||v^{\mu}_{z_1,r}-q_{\mu}||_{\mathcal{L}^2(B^+_1)}+C r\left\|\lambda_-\left( (u_r^{\mu})^{-}\right)^{p-1}-\lambda_+ \left((u_r^{\mu})^{+}\right)^{p-1}\right\|_{\mathcal{L}^{\infty}(B_1')}\right) \leq C(\epsilon).
\end{equation}
for some $C(\epsilon) \to 0$ as $\epsilon \to 0$. We recall that (\ref{Linfty_uni_bnd_1}) and (\ref{Linfty_uni_bnd_2}) hold for all for all $z_1 \in B^+_{\rho_{\epsilon}} (z_0)$ and all $0<r<r_{\epsilon}$. Now, covering the compact set $K$ by finitely many balls of radius $\rho_{\epsilon}$, we can find a small radius $r^K_{\epsilon}$ such that (\ref{Linfty_uni_bnd_1}) and (\ref{Linfty_uni_bnd_2}) hold for all $z_1 \in K$ and for all $0<r<r^K_{\epsilon}$. Clearly, this proves the second part of the theorem.
\end{pf}

In order to study the dimension of the singular set we will need the following result.

\begin{lemma}\label{count_union}
Let $u$ be the solution to (\ref{prb}). Then $\Sigma_{\mu}(u)$ is of type $F_{\sigma}$, that is, it is a union of countably many closed sets.
\end{lemma}
\begin{pf}
Let $$E^{u}_j:=\{ z_0 \in \Sigma_{\mu} (u) \cap \overline{B^+_{1-\frac{1}{j}}} \text { such that } \frac{1}{j} \rho^{\mu} \leq \sup_{B^+_{\rho}(z_0)} |u| <j \rho^{\mu} \}$$
and
$$E^{v}_j:=\{ z_0 \in \Sigma_{\mu} (u) \cap \overline{B^+_{1-\frac{1}{j}}} \text { such that } \frac{1}{j} \rho^{\mu} \leq \sup_{B^+_{\rho}(z_0)} |v| \leq j \rho^{\mu} \}.$$
By Corollary \ref{uv_avg} and Lemma \ref{non_degen} we see that $\Sigma_{\mu}(u)=\left(\cup_{j=1}^{\infty} E^{u}_j\right) \cup \left(\cup_{j=1}^{\infty} E^{v}_j\right)$.
Moreover, in a similar way to  (\cite[Lemma 1.5.3]{GP}), we can show that $E^{u}_j$ and $E^{v}_j$ are closed sets for all $1<j<\infty$.
\end{pf}

Next, we denote by $d^{z_0}_{\mu}$ the dimension of $\Sigma_{\mu}$ at a point $z_0 \in \Sigma_{\mu}$. That is,
\begin{equation}\label{mu_sing_dim}
d^{z_0}_{\mu}=\max \{d^{u,z_0}_{\mu}, d^{v,z_0}_{\mu}\},
\end{equation}
where $$d^{u,z_0}_{\mu}=dim\{\zeta \in \mathbb{R}^{n} \ \ |\ \ \zeta \cdot \nabla_x p^{z_0}_{\mu}(x,0)=0 \text { \ \  for all x } \in \mathbb{R}^n\}$$
and
$$d^{v,z_0}_{\mu}=dim\{\zeta \in \mathbb{R}^{n} \ \ | \ \ \zeta \cdot \nabla_x q^{z_0}_{\mu}(x,0)=0 \text { \ \  for all x } \in \mathbb{R}^n\}.$$
Finally, we let  $\Sigma^d_{\mu}:=\{ z_0 \in \Sigma_{\mu} (u) | d^{z_0}_{\mu}=d \}$. We are now ready to establish the rectifiability of the singular set.
\begin{theorem}\label{singular_dim}
Let $u$ be the solution to (\ref{prb}). Then for every $\mu \in \mathbb{N}$ and $d=0,1,2,..,n-2$, the set $\Sigma^d_{\mu}$ is contained in a countable union of d-dimensional $C^1$-manifold.
\end{theorem}

\begin{pf}
The proof is based on Theorem \ref{cont_dep},  Whitney's extension theorem, and the implicit function theorem. It follows by applying the same arguments as in \cite[Theorem 1.3.8]{GP} to the families of polynomials $\{p_\mu^{z_0}\}$ and  $\{q_\mu^{z_0}\}$ respectively,  to show that $\Sigma^d_{\mu} \cap E^u_j$ and $\Sigma^d_{\mu} \cap E^v_j$ are both contained in a d-dimensional manifold in a neighborhood of $z_0$.
\end{pf}

\textbf{Acknowledgment}
We would like to thank the referees for their careful reading of the manuscript, and for their helpful suggestions.

\end{document}